\documentclass[12pt, reqno]{amsart}

\usepackage{amsmath, amssymb, amsthm, amsfonts}
\usepackage{color}

\newcommand{\F}{\mathcal{F}}
\newcommand{\R}{\mathbb{R}}

\newcommand{\N}{\mathbb{N}}

\newcommand{\weak}{\rightharpoonup}

\newcommand{\half}{\frac{1}{2}}
\newcommand{\norm}[1]{\left\lVert#1\right\rVert}
\newcommand{\abs}[1]{\left\lvert#1\right\rvert}

\newtheorem{thm}{Theorem}[section]
\newtheorem{lemma}[thm]{Lemma}

\newtheorem{defn}[thm]{Definition}

\newtheorem{propn}[thm]{Proposition}
\newtheorem{rem}[thm]{Remark}

\newcommand{\xx}{\mathbf{x}}
\newcommand{\yy}{\mathbf{y}}
\newcommand{\xxx}{\pmb{x}}

\begin{document}
\title[Mass-concentration of blowup solutions to 2D mZK]{Mass-concentration of low-regularity blow-up solutions to the focusing 2D modified Zakharov-Kuznetsov equation}


\author[D. Bhattacharya]{Debdeep Bhattacharya}
\address{
Department of Mathematics \\
Louisiana State University \\
Baton Rouge, LA, USA}
\curraddr{}
\email{debdeepbh@lsu.edu}
\thanks{} 
\date{}
\maketitle 

\begin{abstract}
    We consider the focusing modified Zakharov-Kuznetsov (mZK) equation in two space dimensions. We prove that solutions which blow up in finite time in the $H^1(\R^{2})$ norm have the property that they concentrate a non-trivial portion of their mass (more precisely, at least the amount equal to the mass of the ground state) at blow-up time.
For finite-time blow-up solutions in the $H^s(\R^2)$ norm for $\frac{17}{18} < s < 1$, we prove a slightly weaker result.  Moreover, we prove that the stronger concentration result can be extended to 
the range $ \frac{17}{18} < s \le 1$ 
under an additional assumption on the upper bound of the blow-up rate of the solution.
The main tools used here are the $I$-method and a profile decomposition theorem for a bounded family of $H^1(\R^{2})$ functions.
\end{abstract}

\section{Introduction}%
\label{sec:introduction}

We consider the two-dimensional initial value problem (IVP)
\begin{equation} 
    \begin{cases} v_t + \partial_x (\Delta v)  + \mu \partial_x
	(v^{3}) = 0,\   \xx = (x,y) \in \R^2,\ t >0,\\ v(\xx,0)= v_0(\xx),
    \end{cases} \label{mzk} 
\end{equation} 
where $\mu = \pm 1$, $v$ is a real-valued function, and $\Delta = \partial_x^2 + \partial_y^2$ is the two-dimensional Laplacian operator.
When $\mu=+1$, the equation is known as the modified \textit{focusing} Zakharov-Kuznetsov (mZK) equation in two space dimensions. It is a modification of the standard Zakharov-Kuznetsov (ZK) equation introduced in 3D to model the propagation of nonlinear ion-acoustic waves in magnetized plasma \cite{ZK}. 
On the other hand, when $\mu = -1$, equation \eqref{mzk} is known as the \textit{defocusing} mZK equation.
The mZK equation can also be interpreted as a two-dimensional generalization of the modified Korteweg-de Vries (mKdV) equation, which was deduced by Kakutani and Ono \cite{kaku} to describe the propagation of Alfv\'en waves at a critical angle to the undisturbed magnetic field. 

The mass $\mathcal{M}$ and energy $ \mathcal{E} $ associated with the solution  $v$ to the mZK equation \eqref{mzk} are defined as
\begin{equation*} 
    \mathcal{M}[v(t)] = \int_{\R^2}^{}v^2(\xx,t) d\xx
\end{equation*} 
and
\begin{equation}
    \mathcal{E}[v(t)] = \frac{1}{2}\int_{\R^2} |\nabla v(\xx,t)|^2 d\xx -
    \frac{\mu }{4} \int_{\R^2} v^{4}(\xx,t) d\xx,
    \label{def-energy} 
\end{equation}
respectively. During its lifespan, the solution to the IVP \eqref{mzk} shows conservation of mass and energy. More precisely, any solution $v(t) \in L^2(\R^{2})$ obeys the \textit{mass conservation law}
\begin{equation} \mathcal{M}[v(t)] =
\mathcal{M}[v_0] \label{amass} \end{equation} 
and any solution $v(t) \in H^1(\R^{2})$ obeys the \textit{energy conservation law}
\begin{equation}
\mathcal{E}[v(t)]  =
\mathcal{E}[v_0]
\label{aenergy}
\end{equation}
for all time $t \in [0, T^*)$, where $T^*$ is the maximal time of existence.

The solution to IVP \eqref{mzk} is also invariant under scaling. i.e., if $v$ solves the IVP \eqref{mzk}, then for any $\lambda>0$, the rescaled function $v_\lambda$ defined by
\begin{align}
\label{scaling-gen}
v_{\lambda}(\xx,t) = \lambda v \left(\lambda \xx , \lambda^3 t  \right)
\end{align}	
also solves the equation \eqref{mzk} with initial data
\begin{align*}
v_{0,\lambda}(\xx) = \lambda v_0 \left( \lambda \xx \right).
\end{align*}	

Moreover, note that the $\dot{H}^s(\R^{2})$ norm of the rescaled solution $v_\lambda$ is scaled by $\lambda^{s}$, i.e.
\begin{align*}
\norm{v_\lambda}_{\dot{H}^s(\R^2)} = \lambda^{s} \norm{v}_{\dot{H}^s(\R^2)}.
\end{align*}
 Therefore, the mass of the solution (the $L^2$-norm) is scaling-invariant. For this reason, the 2D mZK equation said to be \textit{mass-critical}.

An interesting property of equation \eqref{mzk} in the focusing case ($\mu = 1$) is that it admits a traveling wave solution in the $x$ direction. 
Let $\varphi$ be a solution to the nonlinear Elliptic equation
\begin{align}
\label{GSeq}
- \Delta_{\xx}  \varphi (\xx) + \varphi (\xx) - \varphi^3 (\xx) = 0,\ \xx \in \R^2.
\end{align}
Then $$ v(\xx,t)=\varphi_c(\xx-ct e_1) $$ is a solution
to the focusing mZK equation which travels only in the $x$-direction with speed $c$, where $e_1$ is the unit vector $(1,0)\in \R^2$, and
$\varphi_c$ is the dilation of $\varphi$ given by $$ \varphi_c(\xx)=\sqrt{c} \,
\varphi(\sqrt{c}\xx) $$ for $c > 0$, and solves the equation $\Delta \varphi_c
-c\varphi_c + \varphi_c^{3} = 0$.
The unique, radial, positive solution $\varphi$ to equation \eqref{GSeq} is known as the \textit{ground state}. The existence of such solution in 2D was shown by
Berestycki, Gallou\"et, and  Kavian \cite{BGK83} (see  Strauss \cite{Sr77},
Berestycki and Lions \cite{BLi83}, and Berestycki, Lions and Peletier \cite{BLP81}
for the existence in other dimensions).
Gidas, Ni, and Nirenberg \cite{gidas-ni-nirenberg} established sufficient conditions to ensure that the positive solutions are radial.
Kwong \cite{Kwong} showed that radial and positive solutions are unique.
The ground state $\varphi$ is appears in the sharp constant of the Gagliardo-Nirenberg inequality (see Weinstein \cite{Weinstein}) 
	\begin{equation}
	    \norm{f}_{L^4(\R^{2})}^4 \le \frac{2}{\norm{\varphi} _{L^2(\R^{2})}^2} \norm{f} _{L^2(\R^{2})}^2 \norm{\nabla f} _{L^2(\R^{2})}^2.
	    \label{GNCrit}
	\end{equation}

The 2D mZK equation has been extensively studied in recent years. 
Biagioni and Linares \cite{BP} studied the local well-posedness in $H^1(\R^2)$.  Linares and Pastor \cite{LPlwp} proved the local well-posedness in $H^s(\R^2)$ for $s>\frac{3}{4}$ and  Ribaud and Vento \cite{RV} improved it to $s>\frac{1}{4}$. More recently, Kinoshita \cite{kinoshita2019well} pushed this result to $s \ge \frac{1}{4}$ and showed that this result is optimal for the Picard iteration approach.

Regarding the global well-posedness, Linares and Pastor \cite{LP}  proved that the focusing 2D mZK equation is globally well-posed in $H^s(\R^2)$ for $s> \frac{53}{63}$ when the mass of the initial data is smaller than the mass of the ground state. Using the $I$-method, Bhattacharya, Farah, and Roudenko \cite{BFR} improved this result to $s > \frac{3}{4}$. Moreover, in the defocusing case, they proved that the same conclusion holds without any assumption on the size of the initial data.

In this paper, our discussion will be limited to the focusing equation only, i.e., when $\mu = +1$.
Unlike in the defocusing case, where all solutions are global in $H^s(\R^2)$ for $\frac{3}{4} <s $, the solutions to the focusing 2D mZK equation are global when the initial data is smaller than the mass of the ground state $ \varphi $ (defined in equation \eqref{GSeq}). However, for initial data with mass bigger or equal to the mass of the ground state, blowup may occur.
Here, we study the behavior of the low-regularity solutions that would exhibit blow up in finite time, which we define as follows:
\begin{defn}
We say that the solution  $u(\xx,t)$ to the IVP \eqref{mzk} with $v_0 \in H^s(\R^{2})$ blows up in finite time if there exists $ 0< T^* < \infty$ such that 
\begin{align*}
\lim_{t \uparrow T^*} \norm{v(\xx,t)}_{H^s(\R^{2})}  = \infty.
\end{align*}	
\end{defn}

Significant amount of work has been done in developing the blowup theory in the last three decades. 
The  2D mZK equation is mass critical, and therefore, can be compared to the critical generalized KdV  (cgKdV) equation in 1D. 
For this equation, Merle \cite{Merle} proved the existence of solutions that blow up in the $H^1(\R^2)$ norm in finite or infinite time. Later, the existence of solutions that blow up in the $H^1(\R^2)$ norm in finite time was established by Martel and Merle \cite{MM}. 
However, far less is known about the blow-up phenomenon for the focusing 2D mZK equation. When the initial  data $v_0$ is in $H^1(\R^{2})$, Farah, Holmer, Roudeko, and Yang \cite{FHRY} proved that there exists $\alpha > 0$ such that the solution to 2D mZK blows up in finite or infinite time if the energy is negative (i.e. $\mathcal{E}[v_0] <0$) and if the mass of the initial data satisfies $\norm{\varphi} _{L^2(\R^{2})} < \norm{v_0}_{L^2(\R^{2})} \le \norm{\varphi} _{L^2(\R^{2})} + \alpha$.  This is referred to as near-threshhold blow-up phenomenon for the negative energy solutions.  Recently, Klein, Roudenko, Stoilov \cite{KRS2020} investigated the $H^1(\R^2)$ blow-up phenomenon for the 2D mZK equation numerically. In particular, they conjectured that for sufficiently localized initial data, with mass larger than the mass of the ground state, the $H^1(\R^2)$ norm of the solution blows up in finite time. In addition, in \cite{KRS2020}, some numerical evidence to support their conjecture is provided.

Kenig, Ponce, and Vega \cite{KPV-con} were the first to establish a concentration phenomenon in the context of the critical gKdV equation. They showed that if a finite-time $H^1(\R)$ blow-up solution exists, there is a universal constant $C_0$ such that the solution concentrates at least $C_0$ amount of its mass at blow-up time.
Mass-concentration results for low-regularity blow-up solutions were first proved by  Colliander, Raynor, Sulem, and Write \cite{CRSW}  for the $L^2$-critical 2D nonlinear Schr\"odinger (NLS) equation. 
Hmidi and Keraani \cite{HK-siam-2006} used a refined version of a compactness lemma adapted to the blowup solution of the $L^2$-critical NLS to improve the results of \cite{CRSW}.
Using a concentration-compactness argument, Tzirakis \cite{tzirakis} proved such results for one dimensional quintic NLS. 
For the cgKdV equation, Pigott \cite{Pigott} proved mass-concentration of $H^1(\R)$ blow-up solutions using the work of Hmidi and Keraani \cite{HK-siam-2006} and the almost conservation law obtained by Farah \cite{FcgKdV}. For the range $ \frac{16}{17} < s < 1$, a slightly weaker result is also proved in \cite{Pigott}, which was strengthened by assuming a precise upper bound on the blow-up rate.

However, regarding the concentration phenomena of blowup solutions to the 2D focusing mZK equation with initial data in the energy space of below (that is, $v_0 \in H^s(\R^2)$ for $s \le 1$), there are no results currently available.
The purpose of this paper is to provide a first analysis of such phenomena.
We prove that a solution that blows up in the $H^s(\R^2)$ norm, for $\frac{17}{18} < s \leq 1$ in finite time has the property that the mass of the solution concentrates inside a ball. 
This property displayed by the finite time blow-up solutions gives more insight into their behavior. 
While our result applies to the $H^1(\R^2)$ blow-up solutions, it is in fact  more general. More precisely, we prove the following result.

\begin{thm}
Let $\frac{17}{18} < s \le 1$. Let $v_0 \in H^s(\R^2)$ and suppose
that the corresponding solution of the focusing 2D mZK equation \eqref{mzk} blows up in finite
time $T^* > 0$.  Let $\gamma(t) > 0$ be such that
\begin{align}
\frac{(T^* - t)^{\frac{s}{3}}}{\gamma(t)} \to 0, 
\label{condition-on-seq-std}
\end{align}
as $ t\uparrow T^*$. Then, there exists $\xx(t) \in \R^2$ such that 
\begin{align}
\limsup \limits_{t \uparrow T^*} \int\limits_{\abs{\xx - \xx(t)} \le \gamma(t)}^{} \abs{v(\xx,t)}^2 d\xx \ge  \norm{\varphi}_{L^2(\R^2)}^2.
\label{limsup-std}
\end{align}

Let $s = 1$,  $\gamma(t)$ is as in \eqref{condition-on-seq-std}, and the same assumptions on $v_0$ and $v(t)$ as above hold. Then there exists $\xx(t) \in \R^2$  such that
\begin{align}
\label{case-s-1}
\liminf \limits_{t \uparrow T^*}  \int\limits_{\abs{\xx  - \xx(t)} \le \gamma(t)}^{} \abs{v(\xx,t)}^2 d\xx \ge  \norm{\varphi}_{L^2(\R^2)}^2.
\end{align}
\label{concentration-cor-std}
\end{thm}

\begin{rem}
    We can choose $\gamma(t) = (T^* - t)^{\frac{s}{3}- \epsilon}$ for any small $\epsilon>0$ satisfying the hypothesis of Theorem \ref{concentration-cor-std}. Thus, the theorem states that as $t$ approaches $T^*$, there exist $\{\xx(t)\}_t\subset \R^2$ such that a non-trivial portion (at least the amount equal to the mass of the ground state) of the mass of the blow-up solution $v(t)$ is concentrated inside a ball centered at $\xx(t)$ with radius $(T^* - t)^{\frac{s}{3}- \epsilon}$,  which shrinks to $0$.
\end{rem}

Note that, in the case $s=1$, Theorem \ref{concentration-cor-std} provides a stronger concentration result, since the limsup in inequality \eqref{limsup-std} is replaced with liminf to obtain \eqref{case-s-1} . We show that such stronger results can also be obtained for $ \frac{17}{18} < s< 1$ if we impose an additional upper bound on the blowup rate of the solution. In particular, we prove the following.

\begin{thm}
Let $\frac{17}{18} < s \le 1$. Let $v_0 \in H^s(\R^{2})$ and supposes that the corresponding solution to IVP \eqref{mzk-sym} with $\mu = +1$ blows up in finite time $T^* > 0$. Suppose, in addition, that there exists $r \ge \frac{1}{3}$ such that 
\begin{align*}
\norm{v(t)} _{H^s(\R^{2})} \lesssim \frac{1}{ \left( T^* - t \right)^{rs}}.
\end{align*}	
Let $\gamma(t)$ satisfy
\begin{align*}
\frac{ \left( T^* - t \right)^{rs}}{\gamma(t)}  \to 0
\end{align*}	
as $t \uparrow T^*$. Then, there exists $\xx(t) \in \R^2$ such that
\begin{align*}
\liminf \limits_{t \uparrow T^*} \int\limits_{\abs{\xx - \xx(t)} \le \gamma(t)}^{} \abs{v(\xx, t)} d\xx \ge \norm{\varphi} _{L^2(\R^{2})}^2.
\end{align*}	
\label{concentration-stronger-std}	
\end{thm}

To prove these results, our strategy is to symmetrize the 2D mZK equation using a linear change of variables.
Following the linear transformation introduced by Gr\"unrock and Herr \cite{GH}, we define the spatial variables
\begin{align}
\label{cov-sym}
\begin{bmatrix}
x' \\ y'
\end{bmatrix} = 
\begin{bmatrix}
a & b \\ a & -b
\end{bmatrix} 
\begin{bmatrix}
x \\ y
\end{bmatrix}
\end{align}
with  $a = 2^{-\frac{2}{3}}$ and $b = 3^{\half}2^{-\frac{2}{3}}$, 
and define $u$ and $u_0$ by 
\begin{align*}
\begin{cases}
u(x',y', t) = v(x,y, t)\\
u_0(x',y') = v_0(x,y).
\end{cases}
\end{align*}	
Then, we have 
\begin{equation}
\begin{cases}
  \partial_x v(x,y) = a (\partial_{x'} + \partial_{y'}) u(x',y')\\
  \partial_y v(x,y) = b (\partial_{x'} - \partial_{y'}) u(x',y'),
\end{cases}
  \label{differentials}
\end{equation}
and
therefore, the focusing 2D mZK equation \eqref{mzk} reduces to the following symmetrized focusing 2D mZK equation
\begin{align}
\begin{cases}
\partial_t u + (\partial_{x'}^3 + \partial_{y'}^3 ) u + a( \partial_{x'} + \partial_{y'}) (u^{3}) = 0,\ \xx' = (x',y') \in \R^2,\ t >0\\
u(\xx',0) = u_0(\xx').
\end{cases}
\label{mzk-sym}
\end{align}

Moreover, using equations \eqref{differentials}, we get
\begin{equation}
  \norm{\nabla v}_{L^2(\R^2)}^2 = \frac{a^2 + b^2}{2ab}
  \norm{\nabla u}_{L^2(\R^2)}^2 + \frac{2(a^2 - b^2)}{2ab} \int_{\R^2}^{} u_{x'}u_{y'} dx'dy'
  \label{covNabla}
\end{equation}
and
\begin{equation}
  \norm{v}_{L^p(\R^2)}^p = \frac{1}{2ab}\norm{u}_{L^p(\R^2)}^p.
  \label{covNabla2}
\end{equation}

Writing the energy $\mathcal{E}[v]$ defined in \eqref{def-energy} in terms of $u$, we get,
\begin{align*}
  \frac{2}{a^2 + b^2} 2ab \mathcal{E}[v] & =  \norm{\nabla u
  }_{L^2(\R^2)}^2 +  \frac{2(a^2 - b^2)}{a^2 + b^2}
  \int_{\R^2}^{} u_{x'}u_{y'} dx' dy' -\frac{2}{4(a^2 +
  b^2)} \norm{u}_{L^4(\R^2)}^4.
\end{align*}
Observing that $a=(a^2 + b^2)^{-1} = 2^{-\frac{2}{3}}$ and $\frac{2(a^2 - b^2)}{a^2 + b^2 } = -1$, and defining the energy of the solution $u$ to the symmetrized focusing 2D mZK equation by 
\begin{equation*}
E[u(t)] = \frac{2ab}{a^2 + b^2} \mathcal{E}[v(t)],
\end{equation*}
we obtain that
\begin{equation}
E[u(t)] = \half \int_{\R^2} |\nabla u(x,y,t)|^2 dxdy - \half \int_{\R^2} (u_xu_y)(x,y,t) dxdy - \frac{ a}{4} \int_{\R^2} u^{4}(x,y,t) dxdy.
  \label{energy}
\end{equation}
Due to \eqref{aenergy}, the energy conservation law for $u(t) \in H^1(\R^2)$ holds, i.e.,
\begin{align*}
    E[u(t)] = E[u_0]
\end{align*}
for all $t \in [0, T^*)$. Moreover, the mass conservation law for $u(t) \in L^2(\R^{2})$ holds due to equation \eqref{amass}.

Once symmetrized, we state the following theorem for the $H^s$-blowup solutions to IVP \eqref{mzk-sym}.  

\begin{thm}
Let $\frac{17}{18} < s \le 1$. Let $u_0 \in H^s(\R^2)$ and suppose
that the corresponding solution of IVP \eqref{mzk-sym} blows up in finite
time $T^* > 0$.  Then, there exists a sequence $t_n \uparrow T^*$ such
that the following statement holds.

There is a function $V \in H^s(\R^2)$ with 
\begin{align}
\label{bound-V}
\norm{V}_{L^2(\R^{2})} \ge \sqrt{2ab} \norm{\varphi}_{L^2(\R^{2})},
\end{align}
a sequence $\{\rho_n\} \in [0, \infty]$, and a family of points $\{\xx'_n\} \subset \R^2$ with 
\begin{align}
\rho_n \le A(T^* - t_n)^{\frac{s}{3}}
\label{rho_n-bound}
\end{align}
for some $A > 0$, such that
\begin{align}
\rho_n u\left( \rho_n \cdot  + \xx'_n ,t_n  \right) \rightharpoonup V \quad
\text{ weakly in } H^s(\R^2) .
\label{weakly}
\end{align}
\label{main}
\end{thm}

As a direct consequence of Theorem \ref{main}, we prove the following concentration result for a $H^s$-blowup solution to the symmetrized equation \eqref{mzk-sym}.
\begin{thm}
Let the hypothesis of Theorem \ref{main} hold. Let $\beta(t) > 0$ 
be such that
\begin{align}
\frac{(T^* - t)^{\frac{s}{3}}}{\beta(t)} \to 0,
\label{condition-on-seq}
\end{align}
as $ t\uparrow T^*$. Then, there exists $\xx'(t) \in \R^2$ such that 
\begin{align}
\limsup \limits_{t \uparrow T^*} \int\limits_{\abs{\xx' - \xx'(t)} \le \beta(t)}^{} \abs{u(\xx',t)}^2 d\xx' \ge 2ab \norm{\varphi}_{L^2(\R^2)}^2.
\label{limsup}
\end{align}

If $s = 1$, there exists $\xx'(t) \in \R^2$ and  $\beta(t)$ as in \eqref{condition-on-seq} such that under that same assumption on $u_0$ and $u(t)$, we have
\begin{align*}
\liminf \limits_{t \uparrow T^*}  \int\limits_{\abs{\xx'  - \xx'(t)} \le \beta(t)}^{} \abs{u(\xx',t)}^2 d\xx' \ge 2ab \norm{\varphi}_{L^2(\R^2)}^2.
\end{align*}
\label{concentration-cor}
\end{thm}
This result is analogous to Theorem \ref{concentration-cor-std}. The only difference is that the mass of the ground state associated with the symmetrized equation \eqref{mzk-sym} is a scalar multiple of that of the un-symmetrized equation \eqref{mzk}. The scalar $2ab$ is associated with Jacobian of the symmetrization transformation \ref{cov-sym}.

Analogous to Theorem \ref{concentration-stronger-std}, for the symmetrized equation, we prove,
\begin{thm}
Let $\frac{17}{18} < s \le 1$. Let $u_0 \in H^s(\R^{2})$ and supposes that the corresponding solution to IVP \eqref{mzk-sym} blows up in finite time $T^* > 0$. Suppose, in addition, that there exists $r \ge \frac{1}{3}$ such that 
\begin{align}
\label{lower-stronger}
\norm{u(t)} _{H^s(\R^{2})} \lesssim \frac{1}{ \left( T^* - t \right)^{rs}}.
\end{align}	
Let $\beta(t)$ satisfies 
\begin{align}
\label{lablablab}
\frac{ \left( T^* - t \right)^{rs}}{\beta(t)}  \to 0
\end{align}	
as $t \uparrow T^*$. 

Then, there exists $\xx'(t) \in \R^2$ such that
\begin{align*}
\liminf \limits_{t \uparrow T^*} \int\limits_{\abs{\xx' - \xx'(t)} \le \beta(t)}^{} \abs{u(\xx', t)} d\xx' \ge 2ab \norm{\varphi} _{L^2(\R^{2})}^2.
\end{align*}	
\label{concentration-stronger}	
\end{thm}

Two of the main ingredients of this work are the $I$-method developed by
Colliander, Keel, Staffilani, Takaoka, and Tao \cite{CKSTT, CKSTTKdV},
and a variant of a compactness theorem in the same spirit as in Hmidi and Keraani \cite{HK-I}.
In Section \ref{sec:notation}, we introduce necessary notations and recall the almost conservation law obtained in Bhattacharya, Farah, and Roudenko \cite{BFR}. In Section \ref{sec:compactness}, we prove a variant of the compactness theorem as a consequence of the usual profile decomposition theorem in $H^1(\R^{2})$.
Section \ref{sec:symm} deals with the symmetrized focusing mZK equation. Following the approach of Hmidi and Keraani \cite{HK-siam-2006} and  Pigott \cite{Pigott}, we prove the concentration results stated in Theorem \ref{main}, Theorem \ref{concentration-cor} and Theorem \ref{concentration-stronger}. Finally, we return to the standard (un-symmetrized) focusing mZK equation in Section \ref{sec:std-mzk} and prove Theorem \ref{concentration-cor-std} and Theorem \ref{concentration-stronger-std}.

\section{Notations and preliminaries}
\label{sec:notation}

In this section, we introduce some notations and recall some of the preliminary results to be used throughout this paper.

The spacial Fourier transform is denoted by $\widehat{(\cdot)}$.  By $\F$, we denote the Fourier transform both in space and time variables. The spacial and temporal frequency variables are denoted by $\zeta=(\xi, \eta)$ and $\tau$, respectively.

By $D^\alpha$ and $J^\alpha$, we define the Fourier multiplier operators with symbols $|\zeta|^\alpha$ and $\langle \zeta \rangle^{\alpha}$, respectively, where $\langle \zeta \rangle = \sqrt{1 + |\zeta|^2}$.
Thus, the norm in the Sobolev space $H^s(\R^2)$ is defined by
\begin{align*}
    \norm{u}_{H^s(\R^2)} = \norm{J^s u}_{L^2(\R^2)}.
\end{align*}

Let $s,\ b \in \R$. The space $X_{s,b}$ is defined as the space of all tempered distributions $u$ on $\R^2\times \R$ such that
\begin{align*}
  \norm{u}_{X_{s,b}} = \norm{ \langle \zeta \rangle^s \langle \tau - \xi^3 - \eta^3 \rangle^b \F(u)(\xi, \eta, \tau)}_{L^2_{\xi,\eta,\tau}} < \infty.
\end{align*}
Also, for $T>0$, we define the restriction norm
\begin{align*}
  \norm{u}_{X_{s,b}^T} = \inf \{ \norm{v}_{X_{s,b}} :  v(t) = u(t) \text{ on } [0,T] \}
\end{align*}

Given $A,B \ge 0$, we write $A \lesssim B$  if for some universal constant $K > 2$, we have $ A \le KB$. We write $A \sim B$ if both $ A \lesssim B$ and $ B \lesssim A$ hold. We write $A << B$ if there is an universal constant $K > 2$ such that $ KA < B$.\\

For arbitrarily small $\varepsilon >0$, we use $a+$ and $a-$ to denote $a+\varepsilon$ and $a-\varepsilon$ respectively. By $a++$ and $a- -$ we denote $a+2\varepsilon$ and $a-2\varepsilon$ respectively.
\\

We define the unitary group associated to the linear part of symmetrized equation (\ref{mzk-sym}) by 
\begin{align*}
  U(t) = e^{-t(\partial_{x'}^3 + \partial_{y'}^3)},
\end{align*}
 i.e., for any $u_0(\xx')$ defined on $\R^2\times \R$, $u(\xx',t) = U(t)u_0(\xx')$ is the solution to the linear IVP
\begin{align*}
  \begin{cases}
  \partial_t u + \partial_{y'}^3 u + \partial_{y'}^3 u = 0\\
  u(\xx',0) = u_0(\xx').
  \end{cases}
\end{align*}

\subsection{The $I$-method}

For a large positive real number $N$ and $0<s<1$, define a Fourier
multiplier operator $I_N: H^s(\R^2)  \rightarrow H^1(\R^2)$ by
\begin{equation*}
  \widehat{I_Nf}(\zeta) = m_N(\zeta) \widehat{f}(\zeta),
\end{equation*}
where $m_N$ is smooth, radially symmetric, non-increasing function of $|\zeta|$ such that
\begin{equation*}
  m_N(\zeta) = 
  \begin{cases}
    1 & \text{ if } |\zeta| \le N\\
    \left( \frac{N}{|\zeta|} \right)^{1-s} & \text{ if } |\zeta| \ge 2N .
  \end{cases}
\end{equation*}

Note that $I_N$ is a smoothing operator of order $(1-s)$ and we have
\begin{align}
\norm{f}_{H^s(\R^2)} \lesssim \norm{I_Nf}_{H^1(\R^2)} \lesssim N^{1-s}
\norm{f}_{H^s(\R^2)}.
\label{deri-gain}
\end{align}

By $\lambda(t)$ and  $\Lambda(t)$, we denote the following quantities related to the $H^s(\R^2)$
norm of the solution $u(\xx',t)$ of IVP \eqref{mzk-sym} 
\begin{align*}
    \lambda(t) & = \norm{u(t)}_{H^s(\R^2)}, \\
    \Lambda(t) & = \sup\limits_{0 \le  \tau \le t} \lambda(\tau).
\end{align*}
Moreover, $\sigma(t)$ and $\Sigma(t)$ denote the quantities related to the $H^1(\R^2)$
norm of the smoothened solution $I_Nu(\xx',t)$.
\begin{align*}
    \sigma(t) & = \norm{I_N u(t)}_{H^1(\R^2)},\\
    \Sigma(t) & = \sup\limits_{0 \le \tau \le t} \sigma(\tau).
\end{align*}
Using, \eqref{deri-gain}, we have the following inequality
\begin{align}
\label{sigma-lambda}
\Sigma(t) \lesssim N^{1-s} \Lambda(t).
\end{align}

Applying the operator $I_N$ on IVP (\ref{mzk-sym}), we get the modified IVP
\begin{equation}\label{igzk}
\begin{cases}
    \partial_t I_Nu + (\partial_{x'}^3 + \partial_{y'}^3) I_Nu +  a (\partial_{x'}+ \partial_{y'}) \left( I_N(u^{3})\right) = 0 \\
  I_Nu(x',y',0) = I_Nu_0(x',y').
\end{cases}
\end{equation}

Now, we prove the following local-wellposedness result for $I_Nu$ similar to Theorem 6.4 of Bhattacharya, Farah, and Roudenko \cite{BFR}. The only improvement here is that the time $\delta$ of local existence of $I_Nu$ depends  on the $-4 - \epsilon$ power of the $H^1(\R^{2})$ norm of the initial data $I_Nu_0$.

\begin{thm}
Let $ \frac{3}{4} < s<1$ and suppose that $u_0 \in H^s(\R^2)$. Then there exists a
$\delta >0$ such that the IVP \eqref{igzk}
has a unique local solution $I_Nu \in C([0, \delta]; H^1(\R^2))$ such
that
\begin{align}
\label{eq:h1bound}
\norm{I_Nu}_{X_{1, \half+}^\delta} \lesssim
\norm{I_Nu_0}_{H^1(\R^2)}.
\end{align}
Moreover, the local time of existence $\delta$ is given by
\begin{align*}
\delta \sim \norm{I_Nu_0}_{H^1(\R^2)}^{-4-}.
\end{align*}
\label{fixed-point}
\end{thm}

\proof
Using the Duhamel's principle we need to find a solution to integral equation
\begin{equation}\label{duhamel}
I_Nu(t) = U(t)I_Nu_0 -  a\int_{0}^{t} U(t-s) (\partial_x+\partial_y) \left( I_N(u^{3}) \right)\, ds.
\end{equation}
To work on $X_{s,b}$ spaces, we consider the following local formulation of Duhamel formula  instead.
\begin{equation}\label{iduhamel}
I_Nu(t) = \psi(t)U(t)I_Nu_0 -  a \psi_\delta(t)\int_{0}^{t} U(t-s) (\partial_x+\partial_y) \left( I_N(u^{3}) \right)\, ds,
\end{equation}
where $\psi \in C_0^\infty([-2,2])$ be an even function with $0 \le \psi \le 1$, $\psi(t) =1$ for $|t| \le 1$,  and $\psi_\delta(t) = \psi(\frac{t}{\delta})$.
Note that if $u: \R^2 \times \R \to \R$ is a solution to equation \eqref{iduhamel}, $u\big|_{\R^2 \times [0,\delta]}$ is a solution to equation \eqref{duhamel} on the interval $[0,\delta]$.

Recalling Lemma 6.1 of \cite{BFR} (See also Lemma 2.2 of \cite{GH}) with $b = \half+$ and $b' = 0$, we have
\begin{equation}\label{linear}
\norm{\psi U(t) u_0}_{X_{s,\half+}} \lesssim \norm{u_0}_{H^s(\R^2)}
\end{equation}
and
\begin{equation}\label{nonlinear}
    \norm{\psi_\delta \int_{0}^{t} U(t-s) f(s) ds }_{X_{s,\half+}} \lesssim \delta^{\half-} \norm{f}_{X_{s,0}}.
\end{equation}
for all $s \in \R$.

Note that in Lemma 3.3 of \cite{BFR}, a stronger version of trilinear estimate is proved, which states that
  for any $\frac{3}{4}<s<1$, 
\begin{equation*}
  \norm{(\partial_x + \partial_y)(u_1 u_2 u_3) }_{X_{s,0}} \le \prod\limits_{i=1}^3 \norm{u_i}_{X_{s,\frac{1}{2}+}}.
  \end{equation*}
Combined with the interpolation lemma stated in 
Lemma 12.1 of \cite{CKSTTKdV}, we obtain
\begin{equation}\label{multilinear}
\norm{(\partial_x + \partial_y) I_N(u ^{3})}_{X_{1,0}}
    \lesssim \norm{I_Nu}_{X_{1,\half +}}^{3}.
\end{equation}

Now, applying $X_{1,\half+}^\delta$ norm on both sides of (\ref{iduhamel})
and applying estimates (\ref{linear}) and (\ref{nonlinear}), we obtain
\begin{equation*}
\norm{I_Nu}_{X_{1,\half+}^\delta} \lesssim
\norm{I_Nu_0}_{H^1(\R^2)} + \delta^{\half-}
\norm{(\partial_x+\partial_y)(I_N(u^{3}))}_{X_{1,0}^\delta}.
\end{equation*}
By definition of localized norm $\norm{\cdot}_{X_{s, b}^\delta}$ , we have 
\begin{equation*}
\norm{I_Nu}_{X_{1,\half+}^\delta} \lesssim
\norm{I_Nu_0}_{H^1(\R^2)} + \delta^{\half-}
\norm{(\partial_x+\partial_y)(I_N(\theta^{3}))}_{X_{1,0}},
\end{equation*}
where the function $\theta(x,y,t) = u(x,y,t)$ on $\R^2 \times [0,\delta]$ and
\begin{equation}
\norm{I_Nu}_{X_{1,0}^\delta} \sim \norm{I_N\theta}_{X_{1,0}}.
\label{equiv-sym}
\end{equation}
Using the estimate (\ref{multilinear}) for $\theta$ and relation (\ref{equiv-sym}), we get
\begin{equation*}
\norm{I_Nu}_{X_{1,\half+}^\delta} \lesssim
\norm{I_Nu_0}_{H^1(\R^2)} +
\delta^{\half-}\norm{I_Nu}_{X_{1,\half+}^\delta}^3.
\end{equation*}

Now, using a standard argument  involving the contraction mapping principle, we conclude that there is $\delta > 0$ with
$\delta^{\half-} \sim \norm{I_N u_0}_{H^1(\R^2)}^{-2}$ such
that  $\norm{I_N u}_{X_{1, \half+}^{\delta}} \lesssim \norm{I_N
u_0}_{H^1(\R^2)}$.
Now the proof is complete.
\qed

The \textit{modified energy} $E^1[u]$ associated with the solution $u$ to the IVP \eqref{mzk-sym} is defined as
\begin{align*}
    E^1[u(t)] = E[I_Nu(t)].
\end{align*}
From Proposition 5.4 of Bhattacharya, Farah, and Roudenko \cite{BFR}, 
we recall the \textit{almost conservation law} for the modified energy.
\begin{propn}
Let $s > \frac{3}{4}$, $N >> 1$ and $u \in H^s(\R^2)$ be a solution to \eqref{mzk-sym} on $[T, T + \delta]$. Then we have the following growth of the modified energy functional
\begin{align}
\label{alm-temp}
\abs{ E^1[u(T+\delta)] - E^1[u(T)] } \lesssim N^{-1+} \left(
\norm{Iu}_{X_{1, \half+}^\delta}^4 + \norm{Iu}_{X_{1,
\half+}^{\delta}}^6 
\right).
\end{align}
\end{propn}

Note that, using \eqref{eq:h1bound} from Theorem \ref{fixed-point}  for $I_Nu$, the almost conservation law \eqref{alm-temp} can be written as
\begin{align}
\label{almostConservation}
\abs{ E^1[u(T+\delta)] - E^1[u(T)] }  \lesssim N^{-1+} \left( \Sigma(T)^4 + \Sigma(T)^6 \right).
\end{align}

\section{A compactness theorem}%
\label{sec:compactness}

In this section, we prove a compactness theorem, which is a variant of Theorem 1.1 of Hmidi-Keraani \cite{HK-siam-2006}. The key difference here is that for a bounded  family of $H^1(\R^{2})$ functions, we have imposed an upper bound on the mass of the transformed gradient under the symmetrization transformation. The trade-off here is that the weak limit of a converging subsequence extracted from this bounded family will be related to the ground-state mass $\norm{\varphi}_{L^2(\R^{2})}$ via a different scalar.
\begin{thm}
Let $\{\psi_n\}$ be a bounded family in $H^1(\R^2)$ such that 
\begin{align}
\label{assump-M}
 \liminf_{n \to \infty} \left(   \norm{\nabla
 \psi_n}_{L^2(\R^2)}^2 - \int\limits_{\R^2}^{} \partial_{x'} \psi_n (\xx') \partial_{y'} \psi_n (\xx') d\xx'  \right) \le M^2
\end{align} 
and
\begin{align}
\label{assump-m}
 \limsup_{n \to \infty} \norm{\psi_n}_{L^4(\R^2)} 
 \ge m.
\end{align} 
Then, there exists $\{\xx'_n \} \subset \R^2$ such that up to a
subsequence, 
\begin{align}
\label{conv-weak}
 \psi_n (.+ \xx'_n) \weak V
\end{align}
weakly in $H^1(\R^2)$ with 
\begin{align}
\label{sharp-bound}
 \norm{V}_{L^2(\R^2)} \ge   \frac{m^2 \sqrt{a^2b}}{M}
 \norm{\varphi}_{L^2(\R^2)},
\end{align}
where $\varphi$ is the unique, positive, radial solution to the elliptic problem
\begin{align*}
\Delta_{\xx'} \varphi(\xx') - \varphi(\xx') + \varphi^3(\xx') = 0,\ \xx' \in \R^2.
\end{align*}
\label{hmidi-keraani}
\end{thm}

The proof of Theorem \ref{hmidi-keraani} relies on the profile decomposition theorem  by Hmidi-Keraani \cite{HK-I} for a bounded sequence in $H^1(\R^n)$, where $n\geq 1$ denotes the spatial dimension. We state the result for $n = 2$. 

\begin{thm}
\label{thm:HK}
Let $ \{ \psi_n \}_{n=1} ^\infty$ be a bounded sequence in $H^1(\R^{2})$. Then, there exist a subsequence, still denoted by $ \{ \psi_n \}_{n=1} ^\infty$, a family $ \{ \xxx'_j \}_{j=1}^\infty$ of sequences in $\R^2$ (where $\xxx'_j = \{ \xx'_{n,j} \}_{n=1}^\infty \subset \R^2$ ), and a sequence $ \{ V^j \}_{j=1} ^\infty$ in $H^1(\R^{2})$, such that the following hold
\begin{enumerate}
\item for every $k \ne j$, $\abs{\xx'_{n,k} - \xx'_{n,j}} \to \infty$ as $n \to \infty$,
\item for every $l \le 1$ and every $\xx' \in \R^2$, we have 
 \begin{align*}
     \psi_n(\xx') = \sum_{j =1}^{l} V^j(\xx' - \xx'_{n,j}) + \psi_n^l(\xx')
 \end{align*}	
	 with
	 \begin{align}
	     \label{eq:l4-conv}
		 \limsup_{n \to \infty} \norm{\psi_n^l} _{L^4(\R^{2})}  \to 0 \text{ as } l \to \infty.
	 \end{align}	
\end{enumerate}
Moreover, the decomposition of both $v_n$ and $\nabla v_n$ satisfy the following asymptotic Pythagorean expansion
\begin{align}
\label{eq:pyth-l2}
 \norm{\psi_n} _{L^2(\R^{2})} ^2 = \sum_{j =1}^{l} \norm{V^j} _{L^2(\R^{2})} ^2 + \norm{\psi_n^l}_{L^2(\R^{2})} ^2 + o(1),
\end{align}
and
\begin{align}
\label{eq:pyth-h1}
 \norm{\nabla \psi_n} _{L^2(\R^{2})} ^2 = \sum_{j =1}^{l} \norm{\nabla V^j} _{L^2(\R^{2})} ^2 + \norm{\nabla \psi_n^l} _{L^2(\R^{2})} ^2 + o(1).
\end{align}
\end{thm}

Here, we prove an additional property for such a decomposition.

\begin{lemma}
In addition to satisfying relations \eqref{eq:pyth-l2} and \eqref{eq:pyth-h1}, the subsequence $ \{ \psi_n \}_{n = 1}^\infty$ and the family  $ \{ V^j \}_{j = 1}^\infty$  in the conclusion of Theorem \ref{thm:HK}  also satisfy the following expansion.
\begin{align}
\label{eq:pyth-mixed}
 \int\limits_{\R^2}^{} \partial_{x'} \psi_n \partial_{y'} \psi_n (\xx') d\xx'  = \sum\limits_{j=1}^{l} \int\limits_{\R^2}^{}  V^j_{x'} (\xx') V^j_{y'} (\xx') d\xx' + \int\limits_{\R^2}^{} \partial_{x'} \psi_n^l (\xx') \partial_{y'} \psi_n^l (\xx') d\xx' + o(1) .
\end{align}
\end{lemma}

\proof
This result is an easy byproduct the proof of Theorem \ref{thm:HK} stated in  Proposition 3.1 of \cite{HK-I}. Here, we present a complete proof by recalling the construction of the families $\{V^j \}_{j = 1}^\infty$ and $ \{ \xxx'_j \}_{j = 1}^\infty $.

Given the bounded sequence $\Psi = \{ \psi_n \}_{n=1} ^\infty \subset H^1(\R^{2})$,  we construct the set $ S(\Psi) $ by collecting all the weak $H^1(\R^{2})$ limits of translated subsequences $\Psi$, i.e.,
\begin{align*}
S(\Psi) = \{ V : \exists \{ \psi_{n_k} \} \subset \Psi \text{ such that }  \psi_{n_k} (\cdot + \xx'_n) \weak V \text{ in } H^1(\R^{2}) \text{ as } k \to \infty , \{ \xx'_n \} \subset \R^2  \}.
\end{align*}	
Next, we take the supremum over the $H^1(\R^{2})$ norm of all such translated subsequencial weak limits. We define
\begin{align*}
\eta(\Psi) = \sup_{V \in S(\Psi)}  \norm{V}_{H^1(\R^{2})} .
\end{align*}	
By definition, we have
$
\eta(\Psi) \le \limsup_{n \to \infty} \norm{\psi_n} _{H^1(\R^{2})} 	.
$

Our aim is to prove that there exists a sequence $ \{ V^j \}_{j=1} ^\infty \subset S(\Psi) $  and a family of sequence of translations $ \{ \xxx'_j \}_{j=1}^\infty $ with a pairwise divergence property
\begin{align*}
\abs{\xx'_{n,k} - \xx'_{n,j}} \to \infty \text{ as } n \to \infty \text{ for } k \ne j 
\end{align*}	
such that upto a subsequence, the sequence $ \{ \psi_n \}_{n=1} ^\infty $ can be decomposed as
\begin{align*}
\psi_n(\xx') = \sum_{j =1}^{l} V^j(\xx' - \xx'_{n,j}) + \psi_n^l(\xx')
\end{align*}	
with
\begin{align*}
\eta(\Psi^l) \to 0 \text{ as } l \to \infty	,
\end{align*}	
where $\Psi^l = \{ \psi_n^l \}_{n=1}^\infty$, 
and the identities \eqref{eq:pyth-l2} - \eqref{eq:pyth-mixed} hold. 
If $\eta(\Psi) = 0$, we an take $V^j = 0$ for all $j$. Otherwise, we choose the first function $V^1 \in S(\Psi)$ such that
\begin{align*}
\norm{V^1} _{H^1(\R^{2})} \ge \half \eta(\Psi) > 0.
\end{align*}	
Since $V^1 \in S(\Psi)$, there exists a sequence of translations $ \{ \xx_{n,1} \}_{n=1} ^\infty$, such that, upto a subsequence, we have
\begin{align*}
\psi_n(\cdot + \xx'_{n,1}) \weak V^1 \text{ in } H^1(\R^{2}).
\end{align*}	
We set  the residue
\begin{align*}
\psi_n^1 : = \psi_n - V^1(\cdot - \xx'_{n,1}).
\end{align*}	
Since $\psi_n^1(\cdot + \xx'_{n,1}) \weak 0$ in $ H^1(\R^{2}) $, as $n \to \infty$,
\begin{align*}
\norm{\psi_n} _{L^2(\R^{2})} ^2 = \norm{V^1} _{L^2(\R^{2})} ^2 + \norm{\psi_n^1} _{L^2(\R^{2})} ^2 + o(1)
\end{align*}	
and 
\begin{align*}
\norm{\nabla \psi_n} _{L^2(\R^{2})} ^2 = \norm{\nabla V^1} _{L^2(\R^{2})} ^2 + \norm{\nabla \psi_n^1} _{L^2(\R^{2})} ^2 + o(1).
\end{align*}	

Moreover, note that
\begin{align*}
&	\int\limits_{\R^2}^{} \partial_{x'} \psi_n (\xx') \partial_{y'} \psi_n (\xx') d\xx'  \\
& =  \int\limits_{\R^2}^{} \partial_{x'}\psi_n(\xx' + \xx'_{n,1}) \partial_{y'} \psi_n ( \xx' + \xx'_{n,1}) d\xx' \\
& = \int\limits_{\R^2}^{}  \partial_{x'}  \left( \psi_n^1(\cdot + \xx'_{n,1}) + V^1 (\xx') \right) \partial_{y'} \left( \psi_n^1 \left( \xx' + \xx'_{n,1} \right) + V^1(\xx') \right) d\xx' \\
& =  \int\limits_{\R^2}^{} \partial_{x'} \psi_n^1(\xx') \partial_{y'} \psi_n^1(\xx') d\xx' + \int\limits_{\R^2}^{} \partial_{x'} V^1(\xx') \partial_{y'} V^1(\xx') d\xx' 
\\ &  	
+ \int\limits_{\R^2}^{} \psi_n^1(\xx' + \xx'_{n, 1}) \partial_{y'} V^1(\xx') d\xx' + \int\limits_{\R^2}^{} \partial_{x'}V^1(\xx') \partial_{y'} \psi_{n,1}(\xx' + \xx'_{n,1}) d\xx'.
\end{align*}	

Recalling the definition of weak convergence, the sequence $ \{ \psi_n \}_{n = 1}^\infty $ weakly converges to $V^1$ in $H^1(\R^{2})$  if and only if  
\begin{align*}
\begin{cases}
\partial_{x'} \psi_n \weak \partial_{x'} V^1 &\text{ in } L^2(\R^{2}), \\
 \partial_{y'} \psi_n \weak \partial_{y'} V^1 &\text{ in } L^2(\R^{2}), \\
 \psi_n \weak \psi  & \text{ in } L^2(\R^{2}).
\end{cases}
\end{align*}	
Therefore, as $n \to \infty$, the integrals  
$
\int\limits_{\R^2}^{} \psi_n^1(\xx' + \xx'_{n, 1}) \partial_{y'} V^1(\xx') d\xx' 
$ 
and
$
\int\limits_{\R^2}^{} \partial_{x'}V^1(\xx') \partial_{y'} \psi_{n,1}(\xx' + \xx'_{n,1}) d\xx'
$
go to zero. Thus, we have
\begin{align*}
\int\limits_{\R^2}^{} \partial_{x'} \psi_n (\xx') \partial_{y'} \psi_n (\xx') d\xx'  
= \int\limits_{\R^2}^{} \partial_{x'} \psi_n^1(\xx') \partial_{y'} \psi_n^1(\xx') d\xx' + \int\limits_{\R^2}^{} \partial_{x'} V^1(\xx') \partial_{y'} V^1(\xx') d\xx'  + o(1).
\end{align*}

Now, replacing $\Psi$ by $\Psi^1$, we repeat the same steps. In other words, if $\eta(\Psi^1) > 0$, we obtain $V^2$  such that its $H^1(\R^{2})$ norm is bigger or equal to half of $\eta(\Psi^1)$. We also get an associated sequence of translations $ \{ \xx'_{n,2} \}_{n=1} ^\infty $ and the set $\Psi^2$. Moreover, we have that
\begin{align}
\label{lim-pf}
\abs{\xx'_{n,1} - \xx'_{n,2}} \to \infty \text{ as } n \to \infty.
\end{align}	
To see  that \eqref{lim-pf} holds, we assume the contrary. i.e., we assume that there exists $\xx_0' \in \R^2$ such that, there is a subsequence of $ \{ \xx'_{n,1} \}_{n = 1}^\infty$ (still denoted by $\xx'_{n,1}$) with
\begin{align*}
\xx'_{n,1} - \xx'_{n,2} \to \xx'_0 \text{ as } n \to \infty.
\end{align*}	
Now, since 
\begin{align*}
\psi_n^1(\xx' + \xx'_{n,2}) = \psi_n^1 \left(\xx' + \left( \xx'_{n,2} - \xx'_{n,1} \right) + \xx'_{n,1}\right),
\end{align*}	
and $\psi_n^1(\cdot + \xx'_{n,1})$ converges weakly to 0, we have that $V^2 = 0$. Thus, $\eta(\Psi^1) = 0$, a contradiction. Therefore, \eqref{lim-pf}  holds.

This way, by iterating the process of orthogonal extraction, we construct the families $ \{ \xxx'_j \}_{j=1} ^\infty $ and $ \{ V^j \}_{j=1}^\infty $ that satisfy the equations \eqref{eq:pyth-l2}, \eqref{eq:pyth-h1}, and \eqref{eq:pyth-mixed}. 

We have chosen $V^j$ such that
\begin{align}
\label{eq:norm-le}
\eta(\Psi^j) \le \norm{V^{j-1}} _{H^1(\R^{2})}
\end{align}	
Since the series $ \sum_{j =1}^{\infty} \norm{V^j}_{H^1(\R^{2})}$ converges,
\begin{align}
\label{eter-le}
\norm{V^j} _{H^1(\R^{2})} \to 0 \text{ as } j \to 0.
\end{align}	
Combining inequalities \eqref{eq:norm-le} and \eqref{eter-le}, we get that $\eta(\Psi^j) \to 0$ as $j \to \infty$.

Therefore, it remains to show that 
\begin{align*}
\limsup_{n \to \infty} \norm{\psi_n^l} _{L^4(\R^{2})}  \to 0 \text{ as } l \to \infty,
\end{align*}	
which follows directly from the proof of Proposition 3.1 of \cite{HK-I}.

Therefore, we have proved the equation \eqref{eq:pyth-mixed} for the subsequence $\{\psi_n\}_{n = 1}^\infty$ and the family $ \{ V^j \}_{j = 1}^\infty $ that satisfy the conclusions of Theorem \ref{thm:HK}. Therefore the proof is complete.
\qed

With the profile decomposition tool developed so far, we are ready to prove 
Theorem \ref{hmidi-keraani}.

\vspace{0.5cm}
\textbf{Proof of Theorem \ref{hmidi-keraani}.}
Using Theorem \ref{thm:HK}, we can find a subsequence of $ \{ \psi_n \}_{n = 1}^\infty $, still denoted by $ \{ \psi_n \}_{n = 1}^\infty $, a family of sequence of translations $ \{ \xx'_{n,j} \}_{n,j=1}^\infty $, a family $\{V^j\}_{n = 1}^\infty$ of $H^1(\R^{2})$ functions such that
\begin{align}
\label{eq:decomp}
\psi_n(\xx') = \sum\limits_{j=1}^{l} V^j(\xx' - \xx'_{n,j}) + \psi_n^l (\xx'),
\end{align}
and \eqref{eq:pyth-l2},  \eqref{eq:pyth-h1} , and \eqref{eq:pyth-mixed}  hold.

Under the change of variables \eqref{cov-sym}, using equations \eqref{covNabla} and \eqref{covNabla2}, the sharp Gagliardo-Nirenberg inequality \eqref{GNCrit} adapted to the variable $\xx'$ is given by
\begin{align}
\label{GN-symm}
\norm{f}_{L^4(\R^2)}^4 \le  \frac{1}{a^2b \norm{\varphi}_{L^2(\R^2)}^2} \left( \norm{\nabla f}_{L^2(\R^2)}^2 - \int\limits_{\R^2}^{}f_{x'}(\xx') f_{y'}(\xx') d\xx' \right) \norm{f}_{L^2(\R^2)}^2.
\end{align}
Taking $f = V^j$ in equation \eqref{GN-symm}, we obtain
\begin{align*}
\norm{V^j}_{L^4(\R^2)}^4 \le C_{sym} \left( \norm{\nabla V^j}_{L^2(\R^2)}^2 - \int\limits_{\R^2}^{} \partial_{x'} V^j (\xx') \partial_{y'} V^j(\xx') d\xx' \right) \norm{V^j}_{L^2(\R^2)}^2,
\end{align*}
where
\begin{align*}
C_{sym} = \frac{1}{ a^2b \norm{\varphi}_{L^2(\R^2)}^2}.
\end{align*}
Summing over $j$, we conclude
\begin{align}
\label{temp-rel}
\sum\limits_{j=1}^{\infty} \norm{V^j}_{L^4\left( \R^2 \right)}^4 \le C_{sym} \left( \sup\limits_{j}  \norm{V^j}_{L^2(\R^2)}^2\right) \sum\limits_{j=1}^{\infty} \left( \norm{\nabla V^j}_{L^2(\R^2)}^2 - \int\limits_{\R^2}^{} \partial_{x'}V^j(\xx') \partial_{y'}V^j(\xx') d\xx' \right).
\end{align}

Subtracting equation \eqref{eq:pyth-mixed} from equation \eqref{eq:pyth-h1}, we obtain
\begin{align}
\label{eq:subt}
\begin{split}
\norm{\nabla \psi_n}_{L^2(\R^2)}^2 - \int\limits_{\R^2}^{} \partial_{x'} \psi_n (\xx') \partial_{y'}\psi_n (\xx') d\xx' = \sum\limits_{j=1}^{l} \left( \norm{\nabla V^j}_{L^2(\R^2)}^2 - \int\limits_{\R^2}^{} \partial_{x'}V^j(\xx') \partial_{y'}V^j(\xx') d\xx' \right) \\
+ \left( \norm{\nabla \psi_n^l}_{2}^2 - \int\limits_{\R^2}^{} \partial_{x'} \psi_n^l (\xx') \partial_{y'} \psi_n^l (\xx') \right) + o(1).
\end{split}
\end{align}
Note that, using Young's inequality, we can write the second term of the right hand side as
\begin{align*}
\norm{\nabla \psi_n^l}_{L^2(\R^{2})}^2 - \int\limits_{}^{} \partial_{x'} \psi_n^l (\xx')  \partial_{y'}\psi_n^l (\xx') d\xx' \sim \norm{\nabla \psi_n^l}_{L^2(\R^{2})}^2 \ge 0.
\end{align*}
Therefore, equation \eqref{eq:subt} and the upper bound \eqref{assump-M} imply
\begin{align*}
& \sum\limits_{j=1}^{\infty} \left( \norm{\nabla V^j}_{L^2(\R^2)}^2 - \int\limits_{\R^2}^{} \partial_{x'}V^j(\xx') \partial_{y'}V^j(\xx') d\xx' \right) \\
& \le \limsup_{n \to \infty} \left( \norm{\nabla \psi_n}_{L^2(\R^2)}^2  - \int\limits_{\R^2}^{} \partial_{x'} \psi_n (\xx') \partial_{y'} \psi_n (\xx') d\xx' \right) 
\le M^2.
\end{align*}
Hence, from \eqref{temp-rel}, we obtain
\begin{align}
\label{eq:M-bound}
\sum\limits_{j=1}^{\infty} \norm{V^j}_{L^4\left( \R^2 \right)}^4 \le C_{sym} M^2 \left( \sup\limits_{j}  \norm{V^j}_{L^2(\R^2)}^2\right).
\end{align}	

On the other hand,  using the triangle inequality, \eqref{eq:decomp}  implies
\begin{align*}
\norm{\psi_n}_{L^4(\R^{2})}^4 \le \norm{ \sum_{j =1}^{l} V^j ( \cdot - \xx'_{n,j}) } _{L^4(\R^{2})} + \norm{\psi_n^l} _{L^4(\R^{2})}.
\end{align*}	
Applying $ \limsup_{n \to \infty} $ on both sides, and using \eqref{eq:l4-conv} and \eqref{assump-m},  we have
\begin{align*}
m^4 \le \limsup_{n \to \infty} \norm{ \sum_{j=1}^{l} V^j(\cdot - \xx'_{n,j})}_{L^4(\R^{2})}^4.
\end{align*}	
Moreover, using the elementary inequality 
\begin{align*}
\abs{ \abs{ \sum_{j =1}^{l} a_j	}^4 - \sum_{j =1}^{l} \abs{a_j}^4 }  \le C \sum_{j \ne k}^{} \abs{a_j} \abs{a_k}^3
\end{align*}	
we can write for all $l \in \N$,
\begin{align*}
& 	 \abs{\norm{ \sum_{j=1}^{l} V^j (\cdot - \xx'_{n,j}) }_{L^4(\R^{2})}^4 - \sum_{j=1}^{l} \norm{V^j (\cdot - \xx'_{n,j})} _{L^4(\R^{2})}} \\
& \le \int\limits_{\R^2}^{}\abs{ \abs{ \sum_{j =1}^{l} V^j (\xx' - \xx'_{n,j}) }^4 - \sum_{j =1}^{l} \abs{V^j (\xx' - \xx'_{n,j}}^4} d\xx' \\
& \le C \sum_{j \ne k}^{} \int\limits_{\R^2}^{}  \abs{V^j(\xx' - \xx'_{n,j})} \abs{V^j \left( \xx' - \xx'_{n,k} \right)}^3 d\xx.
\end{align*}	
Using the pairwise divergence property (conclusion 1 of Theorem \ref{thm:HK}) of the family $ \{ \xxx'_j \} $, we see that the term on the right hand side vanishes as $n \to \infty$.
Therefore, we have the following inequality
\begin{align}
\label{m-bound}
m^4 \le \sum\limits_{j=1}^{\infty} \norm{V^j}_{L^4(\R^{2})}^4.
\end{align}

Combining inequalities \eqref{eq:M-bound} and \eqref{m-bound}, we obtain
\begin{align*}
m^4 \le C_{sym} \left( \sup \limits_{j} \norm{V^j}_2^2 \right) M^2.
\end{align*}
Therefore, we obtain a lower bound on the supremum of the $L^2(\R^{2})$ norm of $V^j$ over all $j \in \N$, i.e.,
\begin{align}
\label{vj-bound}
\sup \limits_{j \in \N} \norm{V^j}_{L^2(\R^{2})} \ge \frac{m^4}{M^2 C_{sym}} .
\end{align}	
Since the series $ \sum_{j=1}^{\infty} \norm{V^j}_{L^2(\R^{2})}^2 $ converges, the supremum is attained at $j=j_0$ (say). Then
\begin{align*}
\norm{V^{j_0}}_{L^2(\R^{2})} \ge \frac{m^2}{\sqrt{C_{sym}} M}. 
\end{align*}	

Changing the variable $\xx'$ to $(\xx' + \xx'_{n,j_0})$, we have
\begin{align*}
\psi_n(\xx' + \xx'_{n,j_0}) = V^{j_0} (\xx') + \sum_{j \ne j_0}^{l} V^j	(\xx' + \xx'_{n,j_0} - \xx'_{n,j}) + \psi_n^l(\xx' + \xx'_{n, j_0}).
\end{align*}	
The pairwise divergence of the family $ \{ \xx'_j \}_{j = 1}^\infty $ implies 
\begin{align*}
V^j(\cdot + \xx'_{n, j_0} - \xx'_{n,j}) \weak 0
\end{align*}	
weakly in $H^1(\R^{2})$ as $n \to \infty$ for every $j \ne j_0$. Therefore, taking the weak $H^1(\R^{2})$ limit, we get
\begin{align*}
\psi_n(\cdot + \xx'_{n,j_0}) \weak V^{j_0} + \psi^l(\cdot + \xx'_{n,j_0}),
\end{align*}	
where $\psi^l$ is the weak $H^1(\R^{2})$ limit of $\psi_n^l$ as $n \to \infty$. However, note that
\begin{align*}
\norm{\psi^l(\cdot + \xx_{n,j_0})} _{L^4(\R^{2})}
\le 
\limsup_{n \to \infty} 	\norm{\psi_n^l(\cdot + \xx_{n,j_0})} _{L^4(\R^{2})}
\le 
\limsup_{n \to \infty} 	\norm{\psi_n^l} _{L^4(\R^{2})}.
\end{align*}	
Since the right hand side goes to zero as $l \to 0$, by uniqueness of weak limits, we obtain
\begin{align*}
\psi^l = 0 
\end{align*}	
for all $ l \ge j_0$. Therefore, 
\begin{align*}
\psi_n(\cdot + \xx'_{n,j_0}) \weak V^{j_0} \text{ in } H^1(\R^{2}).
\end{align*}	
Defining $V:= V^{j_0}$, inequality  \eqref{vj-bound} implies
\begin{align*}
\norm{V}_{L^2(\R^{2})} \ge \frac{m^2}{\sqrt{C_{sym}} M} = \frac{m^2 \sqrt{a^2 b} }{M}  \norm{\varphi}_{L^2(\R^{2})},
\end{align*}
which is the desired bound \eqref{sharp-bound}. Renaming the sequence $\xx'_{n, j_0}$ by $\xx'_n$, \eqref{conv-weak} follows.  
\qed

\section{Mass concentration for blow-up solutions to the symmetrized mZK equation}
\label{sec:symm}

In this section, we prove Theorem \ref{main}, Theorem \ref{concentration-cor}, and Theorem \ref{concentration-stronger} which deal with the symmetrized focusing mZK equation \eqref{mzk-sym}.
But first, we need prove two proposition. The first one provides an upper bound on the modified energy functional $E^1$ in terms of $\Lambda(t)$ and the second one provides a lower bound on the $H^1(\R^{2})$ norm of the smoothened blow-up solution.

\begin{propn}
For $\frac{17}{18}< s \le 1$, there exists $p(s) < 2$ such that the
following statement holds.

Suppose $u_0 \in H^s(\R^2)$ and that $u(t)$ is the corresponding
solution to the IVP \eqref{mzk-sym} on the maximal finite time interval $[0, T^*)$.  Then, for all $T < T^*$ there exists $N = N(T)$
such that the modified energy functional 
\begin{align*}
\abs{E^1_{N(T)}[u(T)]} \le C_0 \Lambda(T)^{p(s)},
\end{align*}
\label{propn:energy}
\end{propn}
where $C_0 = C_0(s, T^*, \norm{u_0}_{H^s})$ and 
\begin{align*}
N(T) = C \left( \Lambda(T) \right)^{\frac{p(s)}{2(1-s)}}.
\end{align*}

\proof
This is a consequence of the almost conservation law \eqref{alm-temp}.
We follow the strategy of Pigott \cite{Pigott}.

First we treat the easy case $s = 1$.
In this case, we take $N(T) = \infty$.
Then, the smoothing operator $I_{N(T)}$ becomes the identity operator, which implies $I_{N(T)} u = u$, and hence, $E^1 = E$. Therefore, the conclusion holds with $p(1) =0$ due to conservation of energy \eqref{energy}  and the fact that $\norm{\nabla u (t)}_{L^2} \to \infty$ and $t \uparrow T^*$.

Next, we suppose that $ \frac{17}{18} < s < 1$. Consider $T$ sufficiently close to $T^*$. Let $N = N (T)$, to be chosen later.
Due to Theorem \ref{fixed-point}, the time of local existence for the solution $u$ to the IVP \eqref{mzk-sym}  is $\delta \sim
\Sigma(T)^{-4 -}$. Therefore, we can divide the interval $[0, T]$ into
$ \frac{T}{\delta}$ sub-intervals of size $\delta$.
On each of those sub-intervals, we use the almost conservation law  \eqref{almostConservation} and obtain
\begin{align*}
E^1[u(T)] & \le E^1[u(0)] + C \frac{T}{\delta} N^{-1+} \left(
 \Sigma(T)^4 + \Sigma(T)^{6}
\right)\\
   & \lesssim N^{2(1-s)} + N^{-1+} \Sigma(T)^{8+} +
   N^{-1+} \Sigma(T)^{10+}.
\end{align*}	
 Here, we have used the Gagliardo-Nirenberg inequality \eqref{GNCrit} to obtain the first term on the right hand side of the inequality above. Indeed, recalling the definition of the energy \eqref{energy}, we get
\begin{align*}
E^1[u(0)] & \lesssim \norm{\nabla I_N u_0}_{L^2(\R^{2})}^2 + \norm{I_N
 u_0}_{L^4(\R^{2})}^4\\
	   & \lesssim  \norm{I_N u_0}_{H^1(\R^{2})}^2 + \norm{I_Nu_0}_{L^2(\R^{2})}^2
	   \norm{\nabla I_N u_0}_{L^2(\R^{2})}^2\\
	   & \lesssim N^{2(1-s)} \norm{u_0}_{H^s(\R^{2})}^2 + N^{2(1-s)} \norm{u_0}_{H^s(\R^{2})}^2 
	   \norm{u_0}_{L^2(\R^{2})}^2\\
	   & \lesssim N^{2(1-s)}.
\end{align*}
 Next, using \eqref{sigma-lambda}, we obtain
 \begin{align}
     \label{eq:test2}
	 E^1[u(T)]  \lesssim N^{2(1-s)} + N^{-1+} N^{8(1-s)+} \Lambda(T)^{8+} +
	 N^{-1+} N^{10(1-s)+} \Lambda(T)^{10+}.
 \end{align}
 Since $N >>1$ and $T$ is close to $T^*$, we have that
 \begin{align*}
     N^{8(1-s)+} \Lambda(T)^{8+} \le N^{10(1-s)+} \Lambda(T)^{10+}.
 \end{align*}	
 Therefore, \eqref{eq:test2} reduces to
 \begin{align*}
     E^1[u(T)]  \lesssim N^{2(1-s)}  +
	 N^{-1+} N^{10(1-s)+} \Lambda(T)^{10+}.
 \end{align*}	

 Now, we choose $N = N(T)$ such that
 \begin{align*}
     N^{2(1-s)} \sim N^{-1+} N^{10(1-s)+} \Lambda(T)^{10+},
 \end{align*}
 or,
 \begin{align}
     \label{choice-N}
	 N(T) \sim \Lambda(T)^{\frac{10}{8s - 7}+}.
 \end{align}
 This choice of $N(T)$ implies
 \begin{align*}
     E^1[u(T)] \lesssim N^{2(1-s)} \sim \Lambda(T)^{p(s)},
 \end{align*}
 where
 \begin{align*}
     p(s) = \frac{20(1-s)}{8s-7}+ .
 \end{align*}
 Notice that $p(s) <2$ if $s > \frac{17}{18}$. Therefore, the result follows.
 \qed

Now, we derive a lower bound on the $H^1(\R^{2})$ norm of the smoothened solution $I_Nu$, where $u$ is a $H^s(\R^{2})$-blowup solution to the symmetrized focusing 2D mZK equation.
\begin{propn}
Let $u_0 \in H^s(\R^2)$ with $0 < s \le 1$. Let $u(t)$ be the solution to IVP \eqref{mzk-sym} on the maximal forward-in-time interval of existence $[0, T^*)$ with $T^* < \infty $. If $t$ is sufficiently close to $T^*$, then 
\begin{align*}
\norm{I_N u(t)}_{H^1(\R^2)} \gtrsim \left( T^* - t \right)^{-\frac{s}{3}}.
\end{align*}
\label{lower-bound}
\end{propn}
\proof
Due to \eqref{deri-gain}, it is enough to show that
\begin{align}
\norm{u(t)}_{H^s(\R^2)} \gtrsim \frac{1}{\left( T^* - t \right)^{\frac{s}{3}}}.
\label{HsLowerBound}
\end{align}

We use a standard scaling argument (see Merle and Raph{\"a}el \cite{Merle-Raphael}, for example).
Note that for the scaled solution defined in equation \eqref{scaling-gen} satisfies
\begin{align*}
\norm{u_\lambda}_{\dot{H}^s(\R^2)} = \lambda^s \norm{u}_{\dot{H}^s(\R^2)}.
\end{align*}
For $t$ sufficiently close to $T^*$, we define $v_t$ by
\begin{align*}
    v_t(\tau, \xx') = \lambda(t) u \left(\lambda(t) x',\ \lambda(t) y',\ t + \lambda^3(t) \tau  \right),
\end{align*}
where $\lambda(t)$ is defined as
\begin{align*}
\lambda(t) = \norm{u(t)}_{\dot{H}^s(\R^{2})}^{-\frac{1}{s} },
\end{align*}	
or,
\begin{align}
    \label{eq:def-lam}
\lambda^s(t) \norm{u(t)}_{\dot{H}^s(\R^2)} = 1.
\end{align}
 Since the equations \eqref{mzk-sym} is $L^2$-critical, equation \eqref{eq:def-lam} implies
\begin{align}
    \label{norm-one}
\norm{v_t(0)}_{\dot{H}^s(\R^2)} = 1.
\end{align}

On the other hand, Theorem 1.1 of \cite{LPlwp} implies that if $u_0 \in H^s(\R^2)$ for $\frac{3}{4} < s \le 1$, there exists $T = T\left( \norm{u_0}_{H^s(\R^2)} \right) > 0$ such that $u(t) \in H^s(\R^2)$ for all $t \in [0, T)$. Note that, even though the result in \cite{LPlwp} concerns the solution $v$ to the standard (un-symmetrized) mZK equation \eqref{mzk}, it applies directly to the solution $v$ to the symmetrized equation \eqref{mzk-sym} since $u$ and $v$ are related via a linear change of variables \eqref{cov-sym}.

Therefore, there exists $\tau_0(s)  > 0$ such that $v_t$ is locally well-posed for all time $\tau \in [0, \tau_0(s))$, which is independent of $t$. This is due to equation \eqref{norm-one} and the fact that the $L^2(\R^{2})$ norm of the initial data $v_t(\xx',0)$ is independent of $t$, i.e.,
\begin{align*}
\norm{v_t(\xx', 0)}_{L^2(\R^{2})} = \norm{u(\xx', 0)}_{L^2(\R^{2})} = \norm{u_0}_{L^2(\R^{2})}.
\end{align*}	
Thus, from the definition of $v_t$,  we conclude that
\begin{align*}
t + \lambda^3(t) \tau_0(s) < T^*.
\end{align*}
Using equation \eqref{eq:def-lam}, we obtain
\begin{align*}
\norm{u(t)}_{\dot{H}^s} \ge \frac{C(s)}{(T^* - t)^{\frac{s}{3}}}
\end{align*}
for some constant $C(s) = \left( \tau_0(s) \right)^{\frac{s}{3} }$.
Therefore, \eqref{HsLowerBound} follows.
\qed

Using the propositions stated above, we proceed to prove Theorem \ref{main}.

 \vspace{0.5cm}
 \textbf{Proof of Theorem \ref{main}.}
 We follow the method described in the proof of Theorem 1.5 of Hmidi and Keraani \cite{HK-siam-2006}. (See also Pigott \cite{Pigott}). Since we consider the 2D symmetrized mZK equation \eqref{mzk-sym}, the new elements here are the definition of $\rho_n$ in terms of the transformed gradient operator, and the use of Theorem \ref{hmidi-keraani}.
 Let $\{t_n\}_{n=1}^\infty$ be a sequence such that $t_n \uparrow T^*$ and for each
 $t_n$ we have $\norm{u(t_n)}_{H^s(\R^2)} = \Lambda(t_n)$.
	 Note that we can choose such a sequence, since $\lambda(t)$ is continuous
	 in $t$, $\Lambda(T^*) = \infty$ and $\Lambda(0) = \norm{u(0)}_{H^s(\R^2)} < \infty$. 


 Define the sequence $u_n$ by 
 \begin{align}
     u_n := \rho_n I_N u ( \rho_n \xx', t_n)
	 \label{vn}
 \end{align}
 and
 \begin{align}
     \rho_n := \frac{\norm{\nabla \varphi}_{L^2(\R^2)}}{\sqrt{\norm{\nabla I_N u(t_n)}_{L^2(\R^2)}^2 - \int\limits_{\R^2}^{} I_N \partial_{x'} u (\xx', t_n) I_N \partial_{y'} u(\xx', t_n)} d\xx'}.
	 \label{rhon}
 \end{align}
 Applying Young's inequality applied to the denominator of $\rho_n$, we note that
 \begin{align*}
     \frac{\norm{\nabla \varphi}_{L^2(\R^{2})}}{\sqrt{\frac{3}{2}} \norm{\nabla I_N u(t_n)}_{L^2(\R^{2})}} \le
	 \rho_n  \le \frac{\norm{\nabla \varphi}_{L^2(\R^{2})}}{\sqrt{\frac{1}{2}} \norm{\nabla I_N u(t_n)}_{L^2(\R^{2})}}.
 \end{align*}

 Since $\norm{u(t_n)}_{H^s(\R^2)}\to \infty$ as $n \to \infty$, using equation \eqref{deri-gain} , we obtain
 $\norm{u(t_n)}_{H^s(\R^{2})} \le \norm{I_Nu(t_n)}_{H^1(\R^2)}$. Therefore, we have
 $\norm{I_Nu(t_n)}_{H^1(\R^2)} \to \infty$ as $n \to \infty$.

 Now, recalling the definition of $H^1(\R^{2})$ norm, 
 \begin{align}
     \label{temp-temp}
	 \norm{I_N u(t_n)}_{H^1(\R^2)} = \norm{I_Nu(t_n)}_{L^2(\R^2)} +
	 \norm{\nabla I_N u(t_n)}_{L^2(\R^2)}.
 \end{align}
 Note that for the first term of the right hand side, we have
 $\norm{I_Nu(t_n)}_{L^2(\R^2)} \le \norm{u_0}_{L^2(\R^2)}$
 due to the definition of the operator $I_N$ and the conservation of mass.
 But the left hand side of \eqref{temp-temp}  approaches infinity as $n \to \infty$. Therefore, we conclude that 
 \begin{align*}\norm{\nabla I_N u(t_n)}_{L^2(\R^2)} \to \infty \text{ as } n \to \infty. \end{align*}
     Therefore,
 \begin{align*}\norm{I_N u(t_n)}_{L^2(\R^2)} \le \norm{\nabla I_N u(t_n)}\end{align*} for some
     $n$ large enough. Thus, for large $n$, we have, 
 \begin{align*}
     \norm{u(t_n)}_{H^s(\R^2)} \lesssim \norm{\nabla I_N u(t_n)}_{L^2(\R^2)},
 \end{align*}
 and hence
 \begin{align}
     \rho_n \lesssim \norm{u(t_n)}_{H^s(\R^2)}^{-1} =
	 \Lambda(t_n)^{-1}.
	 \label{big-lambda}
 \end{align}

 Note that, from Proposition \ref{lower-bound}, there exists a constant $A > 0$ such that $\rho_n \le A(T^* - t_n)^{\frac{s}{3}}$. Therefore, the bound \eqref{rho_n-bound}  holds.


 Due to the scaling invariance of the $L^2$ norm and the fact that $m \le 1$, we observe that
 \begin{align*}
     \norm{u_n}_{L^2(\R^2)} \lesssim \norm{u_0}_{L^2(\R^2)}.
 \end{align*}
 Moreover, from the definition of $u_n$, and using a change of variable, 
 \begin{align*}
     \norm{\nabla u_n}_{L^2(\R^2)}^2 = \rho_n^{4} \int\limits_{\R^2}^{}
	 \abs{ I_N \nabla_{\xx} u (\rho_n \xx', t_n)}^2 d\xx' = \rho_n^2
	 \int\limits_{\R^2}^{} \abs{\nabla I_N u \left( \xx', t_n \right)}^2 d\xx',
 \end{align*}
 and
 \begin{align*}
     \int\limits_{\R^2}^{} \partial_{x'} u_n (\xx') \partial_{y'} u_n
	 \xx' d\xx' 
     & = \rho_n^4 \int\limits_{\R^2}^{}  I_N \partial_{x'} u(
     \rho_n \xx', t_n) (\xx') I_N \partial_{y'} u(\rho_n \xx' ,t_n)  d\xx' \\
     & = \rho_n^2
     \int\limits_{\R^2}^{} I_N \partial_{x'} u(\xx', t_n) I_N \partial_{y'} u(\xx', t_n) d\xx'.
 \end{align*}
 Also, note that
 \begin{align*}
     \norm{u_n}_{L^4(\R^{2})}^4 = \rho_n^4 \int\limits_{\R^2}^{}  \abs{I_Nu(\rho_n
	 \xx', t_n)}^4 d\xx' = \rho_n^2 \int\limits_{\R^2}^{} \abs{I_N u(\xx', t_n)}^4 d\xx'.
 \end{align*}
 So, from Proposition  \ref{propn:energy} and \eqref{big-lambda}, we have that
 \begin{align*}
     E[u_n] = \rho_n^2 E[I_N u(t_n)] \lesssim \rho_n^2 \left(
		 \Lambda(t_n)
	 \right)^{p(s)} \le \left( \Lambda(t_n) \right)^{p(s) -2}.
 \end{align*}
 Therefore,
 $E[u_n] \to 0$ as $n\to \infty$ as 
 $p(s) < 2$ for $\frac{17}{18} < s \le 1$
 and
 $\Lambda(t_n) \to \infty$ as $n \to \infty$.

 From the definition of $\rho_n$,
 \begin{align*}
     & \norm{\nabla u_n}_{L^2(\R^{2})}^2 - \int\limits_{\R^2}^{} \partial_{x'} u_n (\xx') \partial_{y'} u_n (\xx') d\xx' \\
     & = \rho_n^2 \left( \norm{\nabla I_N u(t_n)}_{L^2(\R^{2})}^2  - \int\limits_{\R^2}^{} I_N \partial_{x'} u (\xx', t_n) I_N\partial_{y'} u (\xx', t_n) \right) d\xx' 
     = \norm{\nabla \varphi}_{L^2(\R^{2})}^2.
 \end{align*}
 Therefore, from the definition of $E[u_n]$ and the fact that $E[u_n] \to 0$, we conclude that
 \begin{align*}
     \frac{a}{2} \norm{u_n}_{L^4(\R^{2})}^4 \to \norm{\nabla \varphi}_{L^2(\R^{2})}^2 \text{ as }  n \to \infty.
 \end{align*}

 Now, in Theorem \ref{hmidi-keraani}, we can take
 \begin{align*}
     M = \norm{\nabla \varphi}_{L^2(\R^{2})}
 \end{align*}
 and
 \begin{align*}
     m^2 =  \frac{\sqrt{2}}{\sqrt{a}} \norm{\nabla \varphi}_{L^2(\R^{2})}
 \end{align*}
 to conclude that there exist a function $V \in H^1(\R^{2})$ and a sequence $ \{ \xx'_n \}_{n = 1}^\infty $ such that $u_n (\cdot + \xx'_n) \weak V$   weakly in $H^1(\R^{2})$ with 
 \begin{align*}
     \norm{V}_2 \ge \frac{m^2 \sqrt{a^2b} }{M} \norm{\varphi}_{L^2(\R^{2})} =   \sqrt{2ab} \norm{\varphi}_{L^2(\R^{2})}.
 \end{align*}

 So far, using Theorem  \ref{hmidi-keraani}, we have concluded that there exists a sequence
 $\{\xx'_n\} \subset \R^2$ such that, up to a subsequence, $u_n(\cdot + \xx'_n) \weak V$ weakly in
 $H^1(\R^{2})$ with
 \begin{align*}
     \norm{V}_{L^2(\R^{2})} \ge \sqrt{2ab} \norm{\varphi}_2.
 \end{align*}
 Therefore,  $u_n(\cdot + \xx'_n) \weak V$ weakly in $H^s$ for
 any $s < 1$, which can be rewritten as
 \begin{align*}
     \rho_n I_{N} u( \rho_n \cdot + \xx'_n, t_n) \weak V  
 \end{align*}
 weakly in $H^s$ for any $s < 1$, where we have renamed $\rho_n \xx'_n$
 by $\xx'_n$.


 It remains to show that 
 \begin{align*}
     \rho_n u (\rho_n \cdot + \xx'_n , t_n) \weak V \text{ in } H^s(\R^{2}).
 \end{align*}	

 First, we deal with the easy case when $s =1$.
 When $u_0 \in H^1(\R^{2})$, we choose $N = \infty$. In that case,  $I_{N}  = 1$, the identity operator. Therefore,
 we have the following conclusion:
 \begin{align*} u_n = \rho_n u (\rho_n \cdot + \xx'_n, t_n) \weak V \text{ weakly in }  H^1(\R^2).
 \end{align*}	

 We shall now show that for $ \frac{17}{18} < s < 1$, the $H^s(\R^{2})$ limit of both $\rho_n I_{N} u(\rho_n \cdot + \xx'_n)$ and $\rho_n
 u( \rho_n \cdot + \xx'_n)$ is $V$ (which will imply weak
 convergence of $\rho_n u(\rho_n \cdot + \xx'_n)$ to $V$ in $H^s(\R^{2})$ as well). We do this by bounding the difference of $\rho_n I_N u(\rho_n \cdot + \xx'_n, t_n)$ and $\rho_n u (\rho_n \cdot + \xx'_n, t_n)$ in $\dot{H}^r(\R^{2})$ space with $r < s$.

 Using the scaling property of Fourier transform, and a change of variable,
 \begin{align*}
     \norm{(I_N u - u) (\rho_n \cdot, t_n)}_{\dot{H}^r(\R^{2})} ^2
     & = \frac{1}{\rho_n^4}  \int\limits_{\R^2}^{} \abs{\xi} ^{2r} \abs{\left(m - 1\right) \left( \frac{\xi}{\rho_n}, t_n\right)}^2  \abs{\widehat{u} \left( \frac{\xi}{\rho_n} ,t_n  \right)} ^2 d\xi \\
     & = \rho_n^{2r - 2} \int\limits_{\R^2}^{} \abs{\xi} ^{2r} \abs{ \left( m -1 \right) (\xi)	}  ^2 \abs{ \widehat{u} (\xi, t_n)}^2 d\xi.
 \end{align*}	
 Now, the frequency support of the function $(m -1)$ is $ \{\xi : \abs{\xi} \ge N  \}$. Therefore, for $s>r$,  we can write the integral as
 \begin{align*}
     &
     \rho_n^{2r - 2} \int\limits_{\abs{\xi} \ge N}^{} \abs{\xi} ^{2r} \abs{ \left( m -1 \right) (\xi, t_n)	}  ^2 \abs{ \widehat{u} (\xi, t_n)}^2 d\xi \\
     & =  
 \rho_n^{2r - 2} \int\limits_{\abs{\xi} \ge N}^{} \abs{\xi} ^{2r - 2s} \abs{\xi}^{2s} \abs{ \left( m -1 \right) (\xi, t_n)	}  ^2 \abs{ \widehat{u} (\xi, t_n)}^2 d\xi\\ & \le \rho_n^{2r - 2} N^{2(r-s)} \int\limits_{\abs{\xi} \ge N}^{}  \abs{\xi}^{2s} \abs{ \left( m -1 \right) (\xi, t_n)	}  ^2 \abs{ \widehat{u} (\xi, t_n)}^2 d\xi  \\
     & \le \rho_n^{2r - 2} N^{2(r-s)} \int\limits_{\abs{\xi} \ge N}^{}  \abs{\xi}^{2s}  \abs{ \widehat{u} (\xi, t_n)}^2 d\xi \\
     & \le \rho_n^{2r - 2} N^{2(r-s)} \norm{u(t_n)} _{\dot{H}^s(\R^{2})}^2 .
 \end{align*}	

 Therefore, for any $r < s < 1$, using \eqref{big-lambda} and \eqref{choice-N}, we have
 \begin{align}
     \label{temtem}
	 \norm{\rho_n \left( I_N u - u \right)(\rho_n \cdot +
	 \xx'_n, t_n)}_{\dot{H}^r(\R^2)} \le \rho_n^r N^{r-s} \norm{u(t_n)}_{H^{s}(\R^2)}
	 \lesssim
	 \left( \Lambda(t_n) \right)^{q(s)},
 \end{align}
 where
 \begin{align*}
     q(s) & = \frac{p(s) (r - s)}{2(1-s)} + 1 - r.
 \end{align*}
 Since $\Lambda(t_n) \to \infty$, the right hand side of \eqref{temtem}  goes to zero if
 $q(s) < 0$. Now, $q(s) < 0$ if
 \begin{align*}
     \frac{20(1-s)}{8s - 7} \frac{r -s}{2(1-s)} + 1 -r < 0,
     \text{ or, }
  \frac{10(r -s)}{8s -7} +1 - r < 0.
 \end{align*}
 If $s > \frac{7}{8}$, we have that
 \begin{align}
     \label{r-bound}
     10(r-s) + (8s - 7) (1 - r) < 0, \text{ or, } 
	  r < \frac{2s + 7}{17 - 8s}.
 \end{align}
 Thus,
 if $r < \frac{2s + 7}{17 - 8s}$
 and $s > \frac{7}{8}$,
 \begin{align}
     \label{r-conv}
	 \norm{\rho_n \left( I_N u - u \right)(\rho_n \cdot +
	 \xx'_n, t_n)}_{H^r(\R^2)} \to 0 \text{ as } n \to \infty.
 \end{align}	
 We also assume that $s < \frac{17}{8}$ to
 ensure that the denominator of the upper bound of $r$ in the right hand side of \eqref{r-bound} is positive.
 Further restriction on $s$ from Proposition
 \ref{propn:energy} implies $ s \in (\frac{17}{18} , 1)$. 
 For these values of
 $s$, the range for $r$ becomes $\frac{16}{17} < r < 1$. 
 Therefore,  the convergence
 \eqref{r-conv} holds for any $r \in (\frac{16}{17} , 1)$, in particular, for all $r \in \left(
	 \frac{17}{18}, 1
 \right)$, since $\frac{16}{17} < \frac{17}{18}$. 

 So, $\rho_n I_n u(\rho_n \cdot + \xx'_n, t_n) = \rho_n u (\rho_n \cdot + \xx'_n, t_n) + h_n$ for some $ \{ h_n \}_{n = 1}^\infty \subset H^r(\R^{2})$ with $ r\in (\frac{17}{18} , 1)$ such that $h_n \weak 0$ weakly in $H^r(\R^2)$. This implies 
 $\rho_n u(\rho_n \cdot + \xx'_n, t_n) \weak V$ weakly in $H^s$ for all
 $\frac{17}{18} < s < 1$.

 This completes the proof of Theorem  \ref{main}.
 \qed

 Now, we return to the concentration phenomenon associated with the symmetrized equation \eqref{mzk-sym} , i.e. Theorem \ref{concentration-cor} and Theorem \ref{concentration-stronger}. We follow closely the arguments by Pigott \cite{Pigott}.

 \vspace{0.5cm}
 \textbf{Proof of Theorem \ref{concentration-cor}.}
 First, we assume that $ \frac{17}{18} <  s < 1 $. Using Theorem \ref{main}, we obtain a sequence $\{	t_n\}_{n = 1}^\infty$ with $ t_n \uparrow T^*$ as $n \to \infty$ such that the weak convergence \eqref{weakly} holds. This implies that for all $R > 0$ we have
 \begin{align*}
     \lim\limits_{n \to \infty} \int\limits_{\abs{\xx'} \le R}^{} \rho_n^2 \abs{u ( \rho_n \xx'  + \xx'_n, t_n)}^2 d\xx' \ge \int\limits_{\abs{\xx'} \le R }^{ } \abs{V(\xx')}^2 d\xx'.
 \end{align*}
 Using
 \begin{align*}
     \int\limits_{\abs{\xx'} \le R}^{} \rho_n^2 \abs{u(\rho_n \xx' + \xx'_n, t_n)}^2 d\xx' \le \sup \limits_{\yy \in \R^2} \int\limits_{\abs{\xx'} \le R}^{} \rho_n^2 \abs{u(\rho_n \xx' + \yy, t_n)}^2 d\xx',
 \end{align*}
 applying a change of variable, and taking the limit over $n$, we obtain
 \begin{align}
     \label{ineq-temp}
	 \lim\limits_{n \to \infty}  \sup \limits_{\yy \in \R^2} \int\limits_{\abs{\xx' - \yy} \le \rho_n R}^{} \abs{u( \xx', t_n)}^2 d\xx' \ge \int\limits_{\abs{\xx'}\le R}^{} \abs{V(\xx')}^2 d\xx'.
 \end{align}
 Using the hypothesis that $ \lim\limits_{t \uparrow T^*} \frac{\left( T^* - t \right)^{\frac{s}{3}}}{\beta(t)} = 0$ and the bound \eqref{rho_n-bound} on $\rho_n$, it follows that \begin{align*}
     \lim\limits_{n \to\infty } \frac{\rho_n}{\beta(t_n)} = 0.
 \end{align*}
 Therefore, we can modify the domain of integration of  \eqref{ineq-temp}  to get
 \begin{align*}
     \lim\limits_{n \to \infty}  \sup \limits_{\yy \in \R^2} \int\limits_{\abs{\xx' - \yy} \le \beta(t_n)}^{} \abs{u( \xx', t_n)}^2 d\xx' \ge \int\limits_{\abs{\xx'}\le R}^{} \abs{V(\xx')}^2 d\xx'.
 \end{align*}
 Letting $R \to 	\infty$, and using the bound \eqref{bound-V},  it follows that
 \begin{align}
     \label{last-step}
	 \lim\limits_{n \to \infty}  \sup \limits_{\yy \in \R^2} \int\limits_{\abs{\xx' - \yy} \le \beta(t_n)}^{} \abs{u( \xx', t_n)}^2 d\xx' \ge \int\limits_{\R^2}^{} \abs{V(\xx')}^2 d\xx' \ge 2ab \norm{\varphi}_{L^2(\R^2)}^2.
 \end{align}
 Since $\left\{ t_n \right\}_{n = 1}^\infty$ is a  sequence with $t_n \uparrow T^*$, we have that
 \begin{align*}
     \limsup \limits_{t \uparrow T^*} \sup \limits_{\yy \in \R^2} \int\limits_{\abs{\xx' - \yy} \le \beta(t)}^{}  \abs{u( \xx', t)}^2 d\xx'  \ge 2ab\norm{\varphi}_{L^2(\R^2)}^2.
 \end{align*}
 Using the continuous dependence of the integral $\int\limits_{\abs{\xx' - \yy} \le \beta(t)}^{} \abs{u(\xx', t)}^2  d\xx' $ on $\yy$ for fixed $t \in [0, T)$, and the fact that this integral goes to zero as $\abs{\yy} \to \infty $, we can find $\xx'(t) \in \R^2$  such that
 \begin{align*}
     \sup \limits_{\yy \in \R^2} \int \limits_{\abs{\xx' - \yy}\le \beta(t)}  \abs{u( \xx', t)}^2 d\xx'  = \int\limits_{\abs{\xx' - \xx'(t)}\le \beta(t)}^{} \abs{u(\xx', t)}^2  d\xx'
 \end{align*}
 for all $t \in [0, T)$.
 Therefore \eqref{limsup} follows.

 Next, we prove the result for the case $s=1$. Again, we take $N = \infty$ in the proof of Proposition \ref{lower-bound}, and therefore, $I_N$ becomes the identity  operator. In this case, for any arbitrary sequence $t_n \uparrow T^*$  (as opposed to a particular sequence $t_n$ with $\norm{u(t_n)}_{H^s(\R^2)} = \Lambda(t_n)$ that was used in proof of Theorem \ref{main}) we have the following lower bound on the blowup rate
 \begin{align*}
     \norm{u(t_n)}_{H^1(\R^2)} \gtrsim \left( T^* - t_n \right)^{-\frac{1}{3}}.
 \end{align*}
 In that case, we have the following upper bound on $\rho_n$
 \begin{align*}
     \rho_n \lesssim \left( T^* -t_n \right)^{\frac{1}{3}}.
 \end{align*}
 Using Theorem \ref{main} for $s=1$, we conclude that
 there exists a sequence $\left\{ \xx'_n \right\} \subset \R^2$  and a profile $V \in H^1(\R^2)$ with 
 \begin{align*}
 \norm{V} _{L^2(\R^{2})} \ge \sqrt{2ab} \norm{\varphi}_{L^2(\R^{2})} \end{align*}
 such that for \textit{any} sequence $\left\{ t_n \right\}$ with $t_n \uparrow T^*$, upto a subsequence, we have
 \begin{align*}
     \rho_n u( \rho_n \xx' + \xx'_n, t_n) \rightharpoonup V 
 \end{align*}
 weakly in $H^1(\R^2)$. Since the sequence $\left\{ t_n \right\}$ is arbitrary, following similar steps as in case $s < 1$, we get
 \begin{align*}
     \liminf \limits_{n \to \infty}  \int\limits_{\abs{\xx'} \le R}^{} \rho_n^2 \abs{u( \rho_n \xx' + \xx'_n, t_n)}^2 d\xx' \ge \int\limits_{\abs{\xx'} \le R}^{} \abs{V}^2 d\xx' .
 \end{align*}
 This implies
 \begin{align*}
     \liminf \limits_{t \uparrow T^*} \int\limits_{\abs{\xx' - \xx'(t)} \le \beta(t)}^{} \abs{u(\xx',t)}^2 d\xx' \ge 2ab\norm{\varphi}_{L^2(\R^2)}^2.
 \end{align*}
 Now the proof is complete.
 \qed

 Before we prove Theorem \ref{concentration-stronger}, in the same spirit as in Pigott \cite{Pigott}, we  establish the following proposition to obtain a bound on the modified energy of the blowup solution. It is analogous to Proposition \ref{propn:energy-time}. However, the bound is obtained in terms of $(T^* - T)$ instead of $\Lambda(T)$, when there is an upper bound \eqref{lower-stronger}  on the blowup rate of the solution is imposed. As Theorem \ref{concentration-stronger} suggests,  these blowup solutions concentrate faster compared to the blowup solutions with standard lower bound given by \eqref{HsLowerBound}.
 \begin{propn}
     For $\frac{17}{18}< s \le 1$, there exists $\kappa(s)$ with $\kappa(s)+2rs >0 $ such that the
	 following statement holds.

	 Suppose $u_0 \in H^s(\R^2)$ and that $u(t)$ is the corresponding
	 solution to the IVP \eqref{mzk-sym} on the maximal finite time interval $[0, T^*)$.  Then, for all $T < T^*$ there exists $N = N(T)$
	 such that the modified energy functional satisfies
	 \begin{align}
	     \label{energy-time-bound}
		 \abs{E^1_{N(T)}[u(T)]} \le C_0  (T^* - T)^{\kappa(s)},
	 \end{align}
	 where $C_0 = C_0(s, T^*, \norm{u_0}_{H^s(\R^{2})})$. 
	 Further, $\Lambda(T)$ is related to $N(T)$ by 
	 \begin{align*}
	     N(T) = C(T^* - T)^{\frac{\kappa(s)}{2(1-s)} }.
	 \end{align*}	
	 \label{propn:energy-time}
 \end{propn}

 \proof
 First, using \eqref{deri-gain}, we observe
 \begin{align*}
     \norm{u(t)}_{H^s(\R^{2})} \lesssim N^{1-s} \left( T^* - t \right)^{-rs}.
 \end{align*}	

 We consider $T$ near $T^*$. 
 Due to Theorem \ref{fixed-point}, the time of local existence for the solution $u$ to IVP \eqref{mzk-sym}  is 
 \begin{align*}
     \delta \sim  N^{-4(1-s)-} \left( T^* - T\right)^{rs-}.
 \end{align*}
 Therefore, we can divide the interval $[0, T]$ into
 $ \frac{T}{\delta}$ sub-intervals of size $\delta$ and we use the almost conservation law  \eqref{almostConservation} to get the following increment of the modified energy.
 \begin{align*}
     E^1[u(T)] & \le E^1[u(0)] + C \frac{T}{\delta} N^{-1+} \left(
		 \Sigma(T)^4 + \Sigma(T)^{6}
	 \right)\\
		   & \lesssim N^{2(1-s)} + N^{-1+} N^{8(1-s)+} \left( T^* - T \right)^{-8rs+}
		   + 
		   N^{-1+} N^{10(1-s)+} \left( T^* - T \right)^{-10rs+}.
 \end{align*}	
	 As in Proposition \ref{propn:energy}, have used  the
	 Gagliardo-Nirenberg inequality \eqref{GNCrit}
	 to obtain the first term on the right hand side of the inequality above.

	 Now, we choose $N = N(T)$ such that
	 \begin{align*}
	     N^{2(1-s)} \sim N^{-1+} N^{10(1-s)+} \left( T^* - T \right)^{-10 rs+}.
	 \end{align*}
	 or,
	 \begin{align*}
	     N^{2(1-s) + 1 - 10(1-s) -} \sim \left( T^* -T \right)^{-10rs+},
	 \end{align*}	
	 which implies
	 \begin{align*}
		 N(T) \sim \left( T^* - T \right)^{\frac{-10rs}{8s - 7}+}.
	 \end{align*}
	 This choice of $N(T)$ implies
	 \begin{align*}
	     E^1[u(T)] \lesssim N^{2(1-s)} \sim  \left( T^* - T \right)^{ \frac{-20 rs(1-s)}{8s- 7} +  }.
	 \end{align*}
	 Defining
	 \begin{align*}
	     \kappa(s) = \frac{-20rs(1-s)}{8s-7}+,
	 \end{align*}
	 we notice that \eqref{energy-time-bound}  holds since  $\kappa(s) <2$ if $s > \frac{17}{18}$. Therefore, the result follows.
	 \qed

	 Now we are ready to prove Theorem \ref{concentration-stronger}.

	 \vspace{0.5cm}
	 \textbf{Proof of Theorem \ref{concentration-stronger}}
	 In this proof, we bypass the argument of extracting a maximizing sequence as in the proof of Theorem \ref{main}, thanks to the assumption \eqref{lower-stronger}. Therefore, the conclusion on mass-concentration is stronger than that of \ref{main}.

	 Let $ \{ u_n \}_{n = 1}^\infty $ be a sequence such that $t_n \uparrow T^*$. Therefore, 
	 We define $u_n$ and $\rho_n$ as in the proof of Theorem \ref{main} (see equation \eqref{vn} and \eqref{rhon}).
	 Similar to \eqref{big-lambda}, we can now  bound $\rho_n$ by
	 \begin{align*}
	     \rho_n \lesssim \norm{u(t_n)}_{H^s(\R^{2})}^{-1} \sim (T^* - T)^{rs}.
	 \end{align*}	

	 Further, we have
	 \begin{align*}
	     E[u_n] = \rho_n^2 E[I_N u(t_n)] = \rho_n^2 (T^* - t_n)^{\kappa(s)} \lesssim  \left( T^* - t_n \right)^{\kappa(s) + 2rs}.
	 \end{align*}
	 Since $\kappa(s) + 2rs > 0$, we get $E[u_n] \to 0$ as $n \to \infty$. Arguing as in Theorem \ref{main}, we find a function $V \in H^1(\R^{2})$ and a sequence $ \xx'_n \in \R^2 $	such that
	 \begin{align*}
	     \rho_n u(\rho_n \cdot + \xx'_n) \weak V	 \text{ in } H^s(\R^{2})
	 \end{align*}	
	 and
	 \begin{align*}
	     \norm{V}_{L^2(\R^{2})}^2 \ge 2ab \norm{\varphi}_{L^2(\R^{2})}.
	 \end{align*}	
	 Next, we proceed as in the proof of Theorem \eqref{concentration-cor}.
	 The condition \eqref{lablablab} implies 
	 \begin{align*}
	     \frac{\rho_n}{\beta(t_n)} \to 0 \text{ as } n \to \infty.
	 \end{align*}	
	 Following the proof of Theorem \eqref{concentration-cor} till equation \eqref{last-step}, we have the following statement 
	 \begin{align*}
	     \lim_{n \to \infty} \sup \limits_{\abs{\xx' - \yy} \le \beta(t_n) } \abs{u(\xx', t_n)}^2 d\xx'  \ge 2ab \norm{\varphi}_{L^2(\R^2)}^2.
	 \end{align*}	

	 Since the sequence $\{	t_n \}_{n = 1}^\infty$ is arbitrary, (as opposed to a maximizing sequence), we obtain
	 \begin{align*}
	     \liminf \limits_{t \uparrow T^*} \sup \limits_{\yy \in \R^2} \int\limits_{\abs{\xx' - \yy} \le \beta(t)}^{}  \abs{u( \xx', t)}^2 d\xx'  \ge 2ab\norm{\varphi}_{L^2(\R^2)}^2.
	 \end{align*}
	 This completes the proof of Theorem \ref{concentration-stronger}.
	 \qed

	 \section{Going back to the standard mZK equation}
	 \label{sec:std-mzk}
	 In this section, we prove the concentration results related to the standard (un-symmetrized) focusing 2D mZK equation \eqref{mzk}, i.e., Theorem  
	 \ref{concentration-cor-std} and Theorem \ref{concentration-stronger-std}. These results follow easily by inverting the symmetrization transformation \eqref{cov-sym}.



	 \vspace{0.5cm}
	 \textbf{Proof of Theorem \ref{concentration-cor-std}}
	 By $B_{r}(\xx)$, we denote the closed two dimensional ball of radius $r$ centered at $\xx'$, i.e.
	 \begin{align*}
	     B_{r}(\xx) = \{ \xx' \in \R^2: \abs{\xx' - \xx} \le r\}.
	 \end{align*}	
	 Further, by $B_{a,b}(\xx)$, we denote a closed two dimensional ellipse centered at $\xx$ with the semi-major and semi-minor axes length as $a$ and $b$, respectively. i.e.,
	 \begin{align*}
	     B_{a, b}(\xx) = \left\{ \xx' \in \R^2:  \frac{(x' - x)^2}{a^2} + \frac{(y' - y)^2}{b^2}  \le 1 \right\}.
	 \end{align*}	

	 Recall that inverse transformation of \eqref{cov-sym} is given by 
	 \begin{align}
	     \label{inv-cov}
		 \begin{bmatrix}
		     x \\ y
		 \end{bmatrix} = 
		 \frac{-1}{2ab} 
		 \begin{bmatrix}
		     -b & -b \\ -a & a
		 \end{bmatrix} 
		 \begin{bmatrix}
		     x' \\ y'
		 \end{bmatrix},
	 \end{align}	
	and the Jacobian of the transformation is $2ab$.

	 If $v(t)$ is a solution to \eqref{mzk} with $v_0 \in H^s(\R^{2})$, $ \frac{17}{18} < s < 1$, such that $v(t)$ blows up in finite time $0< T^*< \infty$, defining $ u(\xx', t) := v(\xx,t)$, we see that $u(t)$ is a $H^s(\R^{2})$-blowup solution to the symmetrized 2D mZK \eqref{mzk-sym} with initial data $u_0(\xx') := v_0(\xx)$ and blows up in time $T^*$.

	 Let $\gamma(t)$ be as in \eqref{condition-on-seq-std}. Define
	 \begin{align}
	     \label{gamma-def}
		 \beta(t) := \sqrt{2} a \gamma(t)
	 \end{align}	
	 Note that, since
	 \begin{align*}
	     \frac{(T^* - t)^{\frac{s}{3}}}{\gamma(t)} \to 0 \text{ as } t \uparrow T^*,
	 \end{align*}	
	 we have
	 \begin{align*}
	     \frac{(T^* - t)^{\frac{s}{3}}}{\beta(t)} \to 0 \text{ as } t \uparrow T^*.
	 \end{align*}	

	 Then, Theorem  \ref{concentration-cor} implies that for $ \frac{17}{18} < s < 1$,
	 there exists $\xx'(t) \in \R^2$ such that 
	 \begin{align}
	     \limsup \limits_{t \uparrow T^*} \int\limits_{B_{\beta(t)}(\xx'(t))}^{} \abs{u(\xx',t)}^2 d\xx' \ge 2ab \norm{\varphi}_{L^2(\R^2)}^2.
		 \label{limsup-now}
	 \end{align}
	 From $\xx'(t) = (x'(t), y'(t))$ obtained this way, we define a new sequence $\xx(t) = (x(t), y(t))$ by
	 \begin{align*}
	     \begin{bmatrix}
		 x(t) \\ y(t)
	     \end{bmatrix} = 
		 \frac{-1}{2ab} 
		 \begin{bmatrix}
		     -b & -b \\ -a & a
		 \end{bmatrix} 
		 \begin{bmatrix}
		     x'(t) \\ y'(t)
		 \end{bmatrix}.
	 \end{align*}
	 For a given $\beta(t) > 0$, the relation
	 \begin{align*}
	     \abs{x' - x'(t)}^2  + \abs{y'  - y'(t)}^2 \le \beta^2(t)
	 \end{align*}
	 is equivalent to
	 \begin{align}
	     \label{eq:ellipse}
	     2a^2 \abs{x - x(t)}^2 + 2b^2 \abs{y - y(t)}^2 \le \beta^2(t).
	 \end{align}
The set of points $(x,y)$ satisfying \eqref{eq:ellipse}  can be denoted by the ellipse 
	 $B_{\frac{\beta(t)}{\sqrt{2}a}, \frac{\beta(t)}{\sqrt{2}b}} (\xx(t))$.

	 Using the inverse change of variable \eqref{inv-cov},
	 \begin{align}
	     \label{nownow}
		 \int\limits_{B_{\beta(t)}(\xx'(t))}^{} \abs{u(\xx', t)}^2 d\xx' 
		 = 2ab \int\limits_{B_{\frac{\beta(t)}{\sqrt{2}a}, \frac{\beta(t)}{\sqrt{2}b}} (\xx(t))}^{} \abs{v(\xx, t)}^2 d\xx.
	 \end{align}

	 Substituting \eqref{nownow} in \eqref{limsup-now}, we get
	 \begin{align}
	     \label{result-here}
	     \limsup_{t \uparrow T^*}  \int\limits_{B_{\frac{\beta(t)}{\sqrt{2}a}, \frac{\beta(t)}{\sqrt{2}b}} (\xx(t))}^{} \abs{v(\xx, t)}^2 d\xx \ge \norm{\varphi} _{L^2(\R^{2})}^2.
	 \end{align}	
	 Since $a < b$, using \eqref{gamma-def}, we have
	 \begin{align*}
	     B_{\frac{\beta(t)}{\sqrt{2}a}, \frac{\beta(t)}{ \sqrt{2} b} } \left( \xx (t) \right) \subset B_{\gamma(t)} \left(\xx(t) \right).
	 \end{align*}
	 Therefore, \eqref{result-here} implies
	 \begin{align*}
	     \limsup \limits_{t \uparrow T^*}  \int\limits_{B_{\gamma(t)}\left( \xx'(t)\right)}^{} \abs{v(\xx',t)}^2 d\xx' \ge \norm{\varphi}_{L^2(\R^2)}^2
	 \end{align*}
	 and hence  \eqref{limsup-std}  follows.

	 The case $s=1$ follows similarly.
	 \qed

	 \vspace{0.5cm}
	 \textbf{Proof of Theorem \ref{concentration-stronger-std}.}
	 From Theorem \ref{concentration-stronger}, following the same steps of inverting the symmetrization transformation as in the proof of Theorem \ref{concentration-cor-std},  the result follows.
	 \qed

\subsection*{Acknowledgements}
The author would like to thank Luiz Gustavo Farah for providing important comments and suggestions, which helped improve the manuscript immensely.

		 \bibliography{thesis}

\begin{thebibliography}{10}

\bibitem{BGK83}
H.~Berestycki, T.~Gallou\"{e}t, and O.~Kavian.
\newblock \'{E}quations de champs scalaires euclidiens non lin\'{e}aires dans
  le plan.
\newblock {\em C. R. Acad. Sci. Paris S\'{e}r. I Math.}, 297(5):307--310, 1983.

\bibitem{BLi83}
H.~Berestycki and P.-L. Lions.
\newblock Nonlinear scalar field equations.
\newblock {\em Arch. Rational Mech. Anal.}, 82(4):313--375, 1983.

\bibitem{BLP81}
H.~Berestycki, P.-L. Lions, and L.~A. Peletier.
\newblock An {ODE} approach to the existence of positive solutions for
  semilinear problems in {${\bf R}\sp{N}$}.
\newblock {\em Indiana Univ. Math. J.}, 30(1):141--157, 1981.

\bibitem{BFR}
D.~Bhattacharya, L.~G. Farah, and S.~Roudenko.
\newblock Global well-posedness for low regularity data in the 2d modified
  {Z}akharov-{K}uznetsov equation.
\newblock {\em J. Differential Equations}, 268(12):7962--7997, 2020.

\bibitem{BP}
H.~A. Biagioni and F.~Lgnares.
\newblock Well-posedness results for the modified {Z}akharov-{K}uznetsov
  equation.
\newblock In {\em Nonlinear equations: methods, models and applications
  ({B}ergamo, 2001)}, volume~54 of {\em Progr. Nonlinear Differential Equations
  Appl.}, pages 181--189. Birkh\"{a}user, Basel, 2003.

\bibitem{CRSW}
M.~Chae.
\newblock Ground state mass concentration in the {$L^2$}-critical nonlinear
  {H}artree equation below {$H^1$}.
\newblock {\em Honam Math. J.}, 31(1):45--61, 2009.

\bibitem{CKSTT}
J.~Colliander, M.~Keel, G.~Staffilani, H.~Takaoka, and T.~Tao.
\newblock Almost conservation laws and global rough solutions to a nonlinear
  {S}chr\"{o}dinger equation.
\newblock {\em Math. Res. Lett.}, 9(5-6):659--682, 2002.

\bibitem{CKSTTKdV}
J.~Colliander, M.~Keel, G.~Staffilani, H.~Takaoka, and T.~Tao.
\newblock Multilinear estimates for periodic {K}d{V} equations, and
  applications.
\newblock {\em J. Funct. Anal.}, 211(1):173--218, 2004.

\bibitem{FcgKdV}
L.~G. Farah.
\newblock Global rough solutions to the critical generalized {K}d{V} equation.
\newblock {\em J. Differential Equations}, 249(8):1968--1985, 2010.

\bibitem{FHRY}
L.~G. Farah, J.~Holmer, S.~Roudenko, and K.~Yang.
\newblock Blow-up in finite or infinite time of the 2d cubic
  {Z}akharov-{K}uznetsov equation.
\newblock {\em arXiv preprint arXiv:1810.05121}, 2018.

\bibitem{gidas-ni-nirenberg}
B.~Gidas, W.~M. Ni, and L.~Nirenberg.
\newblock Symmetry of positive solutions of nonlinear elliptic equations in
  {${\bf R}^{n}$}.
\newblock In {\em Mathematical analysis and applications, {P}art {A}}, volume~7
  of {\em Adv. in Math. Suppl. Stud.}, pages 369--402. Academic Press, New
  York-London, 1981.

\bibitem{GH}
A.~Gr\"{u}nrock and S.~Herr.
\newblock The {F}ourier restriction norm method for the {Z}akharov-{K}uznetsov
  equation.
\newblock {\em Discrete Contin. Dyn. Syst.}, 34(5):2061--2068, 2014.

\bibitem{HK-I}
T.~Hmidi and S.~Keraani.
\newblock Blowup theory for the critical nonlinear schr{\"o}dinger equations
  revisited.
\newblock {\em International Mathematics Research Notices},
  2005(46):2815--2828, 2005.

\bibitem{HK-siam-2006}
T.~Hmidi and S.~Keraani.
\newblock Remarks on the blowup for the {$L^2$}-critical nonlinear
  {S}chr\"{o}dinger equations.
\newblock {\em SIAM J. Math. Anal.}, 38(4):1035--1047, 2006.

\bibitem{kaku}
T.~Kakutani and H.~Ono.
\newblock Weak non-linear hydromagnetic waves in a cold collision-free plasma.
\newblock {\em Journal of the physical society of Japan}, 26(5):1305--1318,
  1969.

\bibitem{KPV-con}
C.~E. Kenig, G.~Ponce, and L.~Vega.
\newblock On the concentration of blow up solutions for the generalized {K}d{V}
  equation critical in {$L^2$}.
\newblock In {\em Nonlinear wave equations ({P}rovidence, {RI}, 1998)}, volume
  263 of {\em Contemp. Math.}, pages 131--156. Amer. Math. Soc., Providence,
  RI, 2000.

\bibitem{kinoshita2019well}
S.~Kinoshita.
\newblock Well-posedness for the cauchy problem of the modified
  {Z}akharov-{K}uznetsov equation.
\newblock {\em arXiv preprint arXiv:1911.13265}, 2019.

\bibitem{KRS2020}
C.~Klein, S.~Roudenko, and N.~Stoilov.
\newblock Numerical study of {Z}akharov-{K}uznetsov equations in two
  dimensions.
\newblock {\em arXiv preprint arXiv:2002.07886}, 2020.

\bibitem{Kwong}
M.~K. Kwong.
\newblock Uniqueness of positive solutions of {$\Delta u-u+u^p=0$} in {${\bf
  R}^n$}.
\newblock {\em Arch. Rational Mech. Anal.}, 105(3):243--266, 1989.

\bibitem{LPlwp}
F.~Linares and A.~Pastor.
\newblock Well-posedness for the two-dimensional modified
  {Z}akharov-{K}uznetsov equation.
\newblock {\em SIAM J. Math. Anal.}, 41(4):1323--1339, 2009.

\bibitem{LP}
F.~Linares and A.~Pastor.
\newblock Local and global well-posedness for the 2{D} generalized
  {Z}akharov-{K}uznetsov equation.
\newblock {\em J. Funct. Anal.}, 260(4):1060--1085, 2011.

\bibitem{MM}
Y.~Martel and F.~Merle.
\newblock Blow up in finite time and dynamics of blow up solutions for the
  l2-critical generalized kdv equation.
\newblock {\em Journal of the American Mathematical Society}, 15(3):617--664,
  2002.

\bibitem{Merle}
F.~Merle.
\newblock Existence of blow-up solutions in the energy space for the critical
  generalized kdv equation.
\newblock {\em Journal of the American Mathematical Society}, 14, 07 2001.

\bibitem{Merle-Raphael}
F.~Merle and P.~Raph{\"a}el.
\newblock Blow up of the critical norm for some radial l 2 super critical
  nonlinear schr{\"o}dinger equations.
\newblock {\em American journal of mathematics}, 130(4):945--978, 2008.

\bibitem{Pigott}
B.~Pigott.
\newblock On mass concentration for the critical generalized korteweg--de vries
  equation.
\newblock {\em Proceedings of the Edinburgh Mathematical Society},
  59(2):519--532, 2016.

\bibitem{RV}
F.~Ribaud and S.~Vento.
\newblock A note on the {C}auchy problem for the 2{D} generalized
  {Z}akharov-{K}uznetsov equations.
\newblock {\em C. R. Math. Acad. Sci. Paris}, 350(9-10):499--503, 2012.

\bibitem{Sr77}
W.~A. Strauss.
\newblock Existence of solitary waves in higher dimensions.
\newblock {\em Comm. Math. Phys.}, 55(2):149--162, 1977.

\bibitem{tzirakis}
N.~Tzirakis.
\newblock Mass concentration phenomenon for the quintic nonlinear
  schr{\"o}dinger equation in one dimension.
\newblock {\em SIAM journal on mathematical analysis}, 37(6):1923--1946, 2006.

\bibitem{Weinstein}
M.~I. Weinstein.
\newblock Nonlinear {S}chr\"{o}dinger equations and sharp interpolation
  estimates.
\newblock {\em Comm. Math. Phys.}, 87(4):567--576, 1982/83.

\bibitem{ZK}
V.~E. {Zakharov} and E.~A. {Kuznetsov}.
\newblock {Three-dimensional solitons}.
\newblock {\em Zhurnal Eksperimentalnoi i Teoreticheskoi Fiziki}, 66:594--597,
  Aug. 1974.

\end{thebibliography}
		 \bibliographystyle{abbrv}	

     \end{document}